\documentclass[10pt,reqno]{amsart}
\usepackage{dsfont,amsmath,amsfonts,amscd,amssymb,graphicx,mathrsfs,eufrak,enumitem,euscript,colonequals,upgreek,pifont}
\usepackage[all,arc,curve,color,frame]{xy}
\usepackage[usenames]{color}

\usepackage{setspace}
\usepackage{array,multirow,booktabs,longtable}

\setlist[enumerate]{format=\normalfont}

\usepackage[colorlinks]{hyperref}

\hypersetup{
  colorlinks   = true,          
  urlcolor     = blue,          
  linkcolor    = purple,          
  citecolor   = blue             
}
\usepackage{tikz,mathrsfs}
\usetikzlibrary{arrows,decorations.pathmorphing,decorations.pathreplacing,positioning,shapes.geometric,shapes.misc,decorations.markings,decorations.fractals,calc,patterns}

\tikzset{
        cvertex/.style={circle,draw=black,inner sep=1pt,outer sep=3pt},
        vertex/.style={circle,fill=black,inner sep=1pt,outer sep=3pt},
        vertexblank/.style={circle,fill=black!40!white,inner sep=1pt,outer sep=3pt},
        DBs/.style={circle,draw=black,circle,fill=green!50!black,inner sep=0pt, minimum size=3pt},
         DWs/.style={circle,draw=black,circle,fill=white,inner sep=0pt, minimum size=3pt},
         DWds/.style={circle,draw=black,densely dotted,circle,fill=white,inner sep=0pt, minimum size=3pt},
        tvertex/.style={inner sep=1pt,font=\scriptsize},
        gap/.style={inner sep=0.5pt,fill=white},
        Ggap/.style={inner sep=0.5pt,fill=green!40!black!20}}

\tikzset{
        DB/.style={circle,draw=black,circle,fill=white,inner sep=0pt, minimum size=5pt},
        DW/.style={circle,draw=red,fill=red,inner sep=0pt, minimum size=5pt},
        DWr/.style={circle,draw=blue,fill=blue,inner sep=0pt, minimum size=5pt},
        dotted/.style={circle,draw=black,densely dotted,fill=white,inner sep=0pt, minimum size=5pt},
}

\tikzstyle{mybox} = [draw=black, fill=blue!10, very thick,
    rectangle, rounded corners, inner sep=10pt, inner ysep=20pt]
\tikzstyle{boxtitle} =[fill=blue!50, text=white,rectangle,rounded corners]

\newcommand{\arrow}[2][20]
 {
  \hspace{-5pt}
  \begin{tikzpicture}
   \node (A) at (0,0) {};
   \node (B) at (#1pt,0) {};
   \draw [#2] (A) -- (B);
  \end{tikzpicture}
  \hspace{-5pt}
 }

\newcommand{\birational}[1][20]{\arrow[#1]{->,dashed}}

\newtheorem{thm}{Theorem}[section]

\newtheorem{defin}[thm]{Definition}
\newtheorem{cor}[thm]{Corollary}

\newtheorem{conj}[thm]{Conjecture}

\theoremstyle{definition} 

\newtheorem{example}[thm]{Example}

\newtheorem{remark}[thm]{Remark}

\newtheorem{notation}[thm]{Notation}

\numberwithin{equation}{section}

\newcommand\citetype[1]{}

\newcommand{\m}{\mathfrak{m}}
\newcommand{\n}{\mathfrak{n}}

\def\Cl{\mathop{\rm Cl}\nolimits}

\def\CM{\mathop{\rm CM}\nolimits}

\def\uCM{\mathop{\underline{\rm CM}}\nolimits}

\def\depth{\mathop{\rm depth}\nolimits}

\def\mod{\mathop{\rm mod}\nolimits}

\def\coh{\mathop{\rm coh}\nolimits}
\def\Qcoh{\mathop{\rm Qcoh}\nolimits}
\def\Mod{\mathop{\rm Mod}\nolimits}

\def\refl{\mathop{\rm ref}\nolimits}

\def\pd{\mathop{\rm pd}\nolimits}

\def\uEnd{\mathop{\underline{\rm End}}\nolimits}
\def\Hom{\mathop{\rm Hom}\nolimits}

\def\RHom{\mathop{\rm {\bf R}Hom}\nolimits}
\def\End{\mathop{\rm End}\nolimits}
\def\Ext{\mathop{\rm Ext}\nolimits}

\def\add{\mathop{\rm add}\nolimits}

\def\Supp{\mathop{\rm Supp}\nolimits}

\def\Spec{\mathop{\rm Spec}\nolimits}

\def\gl{\mathop{\rm gl.dim}\nolimits}

\def\Db{\mathop{\rm{D}^b}\nolimits}

\def\Id{\mathop{\rm{Id}}\nolimits}

\newcommand{\K}{\mathop{{}_{}\mathds{k}}\nolimits}

\newcommand{\con}{\mathrm{con}}
\newcommand{\aff}{\mathrm{aff}}

\newcommand{\CA}{\mathrm{A}_{\con}}

\def\ab{\mathop{\rm ab}\nolimits}

\def\uotimes{\mathop{\underline{\otimes}}\nolimits}

\def\RHom{{\rm{\bf R}Hom}}

\def\sEnd{\scrE\mathrm{nd}}

\newcommand\RDerived[1]{{\rm\bf R}{#1}}

\newcommand\art{\mathsf{Art}}

\newcommand\cart{\mathsf{CArt}}

\newcommand\Sets{\mathsf{Sets}}
\newcommand\cDef{\mathrm{c}\scrD\mathrm{ef}}
\newcommand\Def{\scrD\mathrm{ef}}

\newcommand\DG{\mathsf{DG}}

\newcommand{\cX}{\mathcal{X}}
\newcommand{\cY}{\mathcal{Y}}

\newcommand\twistGen{{\sf Twist}}

\newcommand{\scrA}{\EuScript{A}}

\newcommand{\scrC}{\EuScript{C}}
\newcommand{\scrD}{\EuScript{D}}
\newcommand{\scrE}{\EuScript{E}}
\newcommand{\scrF}{\EuScript{F}}
\newcommand{\scrG}{\EuScript{G}}
\newcommand{\scrH}{\EuScript{H}}
\newcommand{\scrI}{\EuScript{I}}
\newcommand{\scrJ}{\EuScript{J}}
\newcommand{\scrK}{\EuScript{K}}
\newcommand{\scrL}{\EuScript{L}}

\newcommand{\scrN}{\EuScript{N}}
\newcommand{\scrO}{\EuScript{O}}

\newcommand{\scrR}{\EuScript{R}}

\newcommand{\scrW}{\EuScript{W}}
\newcommand{\scrX}{\EuScript{X}}
\newcommand{\scrY}{\EuScript{Y}}
\newcommand{\scrZ}{\EuScript{Z}}

\newcommand{\Jac}{\scrJ\mathrm{ac}}
\def\GKdim{\mathop{\rm JRdim}\nolimits}
\newcommand{\WDelt}{W_{\kern -0.1em \Updelta}}
\newcommand{\Cone}[1]{{\sf Cone}{ (#1)}}
\newcommand{\cham}{\operatorname{\sf Cham}\nolimits}
\newcommand{\Wkern}[1]{W_{\kern -0.1em #1}\kern 0.05em}
\newcommand{\Level}{\operatorname{\sf Level}}
\newcommand{\Delt}{\Updelta}

\newcommand{\sigE}{{^{\upsigma}\kern -1pt\scrE}}
\newcommand{\nsigE}{{^{{\scriptsize-}\upsigma}\kern -1pt\scrE}}

\def\Id{\mathop{\rm{Id}}\nolimits}

\def\Db{\mathop{\rm{D}^b}\nolimits}

\def\add{\mathop{\rm add}\nolimits}

\numberwithin{equation}{section}

\newcommand{\llangle}{\langle\!\langle}
\newcommand{\rrangle}{\rangle\!\rangle}
\newcommand{\llsq}{[\![}
\newcommand{\rrsq}{]\!]}
\newcommand{\llbr}{(\!(}
\newcommand{\rrbr}{)\!)}

\newcommand{\Acon}{\mathrm{A}_{\textnormal{con}}}
\newcommand{\Bcon}{\mathrm{B}_{\textnormal{con}}}

\makeatletter
\newcommand*\bigcdot{\mathpalette\bigcdot@{.5}}
\newcommand*\bigcdot@[2]{\mathbin{\vcenter{\hbox{\scalebox{#2}{$\m@th#1\bullet$}}}}}
\makeatother

\newcommand{\Esix}[6]{%
\begin{tikzpicture}[scale=0.21]
\node at (0,0) [#1] {};
\node at (1,0) [#2] {};
\node at (2,0) [#3] {};
\node at (2,1) [#4] {};
\node at (3,0) [#5] {};
\node at (4,0) [#6] {};
\end{tikzpicture}
}

\newcommand{\Eeight}[8]{%
\begin{tikzpicture}[scale=0.21]
\node at (0,0) [#1] {};
\node at (1,0) [#2] {};
\node at (2,0) [#3] {};
\node at (2,1) [#4] {};
\node at (3,0) [#5] {};
\node at (4,0) [#6] {};
\node at (5,0) [#7] {};
\node at (6,0) [#8] {};
\end{tikzpicture}
}
\tikzset{
W/.style={circle,draw=black,circle,fill=white,inner sep=0pt, minimum size=4pt},
B/.style={circle,draw=black!80!white,circle,fill=black!80!white,inner sep=0pt, minimum size=4pt},
R/.style={circle,draw=black!80!white,circle,fill=red!80!white,inner sep=0pt, minimum size=4pt},  
}

\begin{document}

\title{\textsc{A lockdown survey on cDV singularities}}
\author{Michael Wemyss}
\address{Michael Wemyss, The Mathematics and Statistics Building, University of Glasgow, University Place, Glasgow, G12 8QQ, UK.}
\email{michael.wemyss@glasgow.ac.uk}
\begin{abstract}
This is an expository survey article on compound Du Val (cDV) singularities, with emphasis on recent homological approaches, including: noncommutative resolutions, tilting theory, contraction algebras, classification, derived categories, autoequivalences, stability conditions, deformation theory, and cluster tilting theory.
\end{abstract}
\thanks{This survey forms part of the authors 2019/20 Adams Prize, from the University of Cambridge}
\maketitle
\parindent 20pt
\parskip 0pt

\section{Introduction}\label{L0}

Introduced by Reid \cite{Pagoda} in the early 1980s, compound Du Val (=cDV) singularities can be viewed as a $3$-dimensional `lifting' of the more famous $2$-dimensional Du Val, or equivalently Kleinian, singularities.  It is fair to say, from many algebraic and geometric viewpoints, that cDV singularities are equally as fundamental.  The purpose of this survey is to give an overview of their equivalent definitions and properties, with emphasis on recent homological and noncommutative approaches.

\medskip
In dimension two, Kleinian singularities are key objects in algebraic geometry.  Much of their power comes from their ubiquity: they have strong links with Lie theory, symplectic singularities, representation theory, groups, quivers, and many other categorical and noncommutative structures.  In comparison, cDV singularities have been much less extensively studied, outside of the context of the minimal model programme (MMP).   Of course, reasons abound: there is no good and easy `list' of cDV singularities, there are uncountably many of them, they are neither graded nor always isolated, their representation theory is almost always wild, and their birational geometry is more complicated than that for surfaces.  Basically, they are just harder. However, all these problems can be overcome, and the aim of this survey is to explain how.

\subsection{Definition via the Elephant}
The first, and most elementary, definition of cDV singularities is the following explicit analytic characterisation.
\begin{defin}\label{cDV def intro}
A complete local $\mathbb{C}$-algebra $\scrR$ is called a \emph{cDV} singularity if
\[ 
\scrR\cong\frac{\mathbb{C}\llsq u,v,x,t\rrsq}{f+tg}
\]
where $f\in\mathbb{C}\llsq u,v,x\rrsq$ defines a Du Val, or equivalently Kleinian, surface singularity and $g\in\mathbb{C}\llsq u,v,x,t \rrsq$ is arbitrary.
\end{defin}

In particular, cDV singularities are complete locally (in fact, stalk locally) hypersurfaces, and clearly factoring by $t$ recovers the Kleinian singularity $f$.  This factoring, to recover a Kleinian singularity, is a particular example of the \emph{elephant} \cite{Pagoda}.

More globally, we say that a variety has at worst cDV singularities if, at every closed point, the completion of the stalk satisfies \ref{cDV def intro}.  This definition, whilst being local and explicit, has clear drawbacks.  The isomorphism in the displayed equation forces us to work up to formal changes in coordinates, and recognising equations up to this is both challenging and unpleasant.  Further, knowing an equation is not usually the first thing we want to know about a variety; we are usually more interested in its properties.  Thankfully the abstraction of birational geometry largely fixes this, and is explained and expanded upon below.  However, the definition does serve to illustrate two clear points, both of which are a consequence of the arbitrary choice of $g$, namely (1) cDV singularities, unlike their surfaces counterparts, have moduli, and (2) case-by-case checks will not cut the mustard here.  There are, however, still some special cases that give a reasonable intuition for the general case. One such family is explained in \S\ref{Type A section}.

There are some normal forms for the equations of cDV singularities provided by Markushevich \cite{Markushevich}, however these often turn out to have limited value, given most of the information of a cDV singularity takes place in its birational geometry, and this cannot be read directly from its equations.  Indeed, the best way to understand cDV singularities is through their role in the MMP, which is where they were first discovered.

\subsection{Geometric Overview}\label{geo overview intro}
We outline some of the key features of the birational geometry here, delaying some of the more technical aspects, such as the characterisation of Gorenstein terminal singularities, until \S\ref{geometry section}.   

Kleinian singularities are surfaces, and as such admit a minimal (or, crepant) resolution.  The fibre above the origin is a tree of curves, in ADE formation.  Moving up to $3$-folds, the world of cDV singularities is significantly more complicated.  In brief: (1) crepant resolutions are no longer the aim, as the `best answer' are \emph{minimal models}, which are often singular, (2) these minimal models are no longer unique, as they are related by surgeries called flops, and (3) the fibre above the origin is still a collection of curves, but the dual graph is neither ADE nor even a tree in general, and the formation of the curves can vary between the different minimal models.  Still, the upshot in crude terms is that we input a cDV singularity $\scrR$, and output a finite number of minimal models $\scrX_i\to\Spec\scrR$. 

All minimal models of cDV singularities (which always exist), and more generally all crepant partial resolutions have fibres of dimension at most one.  Hence necessarily $\scrX\to\Spec\scrR$ is one of two types: either a curve-to-point, or a divisor-to-curve contraction.  As a cartoon, the simplest case is when there is a single curve above the origin, and the difference between the two cases is roughly sketched below:
\begin{equation}
\begin{array}{ccc}
\begin{array}{c}
\begin{tikzpicture}
\node at (-1.7,0) {$\scrX$};
\node at (-1.9,-2) {$\Spec\scrR$};
\node at (0,0) {\begin{tikzpicture}[scale=0.5]
\coordinate (T) at (1.9,2);
\coordinate (B) at (2.1,1);
\draw[red,line width=1pt] (T) to [bend left=25] (B);
\draw[color=blue!60!black,rounded corners=5pt,line width=1pt] (0.5,0) -- (1.5,0.3)-- (3.6,0) -- (4.3,1.5)-- (4,3.2) -- (2.5,2.7) -- (0.2,3) -- (-0.2,2)-- cycle;
\end{tikzpicture}};
\node at (0,-2) {\begin{tikzpicture}[scale=0.5]
\filldraw [red] (2.1,0.75) circle (1pt);
\draw[color=blue!60!black,rounded corners=5pt,line width=1pt] (0.5,0) -- (1.5,0.15)-- (3.6,0) -- (4.3,0.75)-- (4,1.6) -- (2.5,1.35) -- (0.2,1.5) -- (-0.2,1)-- cycle;
\end{tikzpicture}};
\draw[->, color=blue!60!black] (0,-1) -- (0,-1.5);
\node at (0,-2.75) {\ding{192}};
\end{tikzpicture}
\end{array}
& \mbox{or} &
\begin{array}{c}
\begin{tikzpicture}
\node at (0,0) {\begin{tikzpicture}[scale=0.5]
\coordinate (T) at (1.9,2);
\coordinate (B) at (2.1,1);
\draw[red,line width=1pt] (T) to [bend left=25] (B);
\foreach \y in {0.1,0.2,...,1}{ 
\draw[very thin,red!50] ($(T)+(\y,0)$) to [bend left=25] ($(B)+(\y,0)$);
\draw[very thin,red!50] ($(T)+(-\y,0)$) to [bend left=25] ($(B)+(-\y,0)$);;}
\draw[color=blue!60!black,rounded corners=5pt,line width=1pt] (0.5,0) -- (1.5,0.3)-- (3.6,0) -- (4.3,1.5)-- (4,3.2) -- (2.5,2.7) -- (0.2,3) -- (-0.2,2)-- cycle;
\end{tikzpicture}};
\node at (0,-2) {\begin{tikzpicture}[scale=0.5]
\draw [red] (1.1,0.75) -- (3.1,0.75);
\draw[color=blue!60!black,rounded corners=5pt,line width=1pt] (0.5,0) -- (1.5,0.15)-- (3.6,0) -- (4.3,0.75)-- (4,1.6) -- (2.5,1.35) -- (0.2,1.5) -- (-0.2,1)-- cycle;
\end{tikzpicture}};
\draw[->, color=blue!60!black] (0,-1) -- (0,-1.5);
\node at (0,-2.75) {\ding{193}};
\end{tikzpicture}
\end{array}
\end{array}\label{two cases intro}
\end{equation}
In general, the fibres can have an arbitrary number of components, and $\scrX$ need not be smooth.  Some other pictures, which are also very misleading but for wildly different reasons, are provided in \S\ref{min models sect}, \S\ref{Type A section} and \cite[5.5]{IW5}.

The combinatorics of the birational geometry of cDV singularities is controlled not by Lie theory, but by \emph{intersection arrangements} inside ADE root systems, and we explain this in \S\ref{combo section}.  These hyperplane arrangements need not even be Coxeter in general, and so the language describing them becomes more clunky, and less elegant.  An example of such an intersection arrangement is the following, with 192 chambers in $\mathbb{R}^3$.  
\begin{equation}
\begin{array}{c}
\includegraphics[angle=0,scale = 0.2]{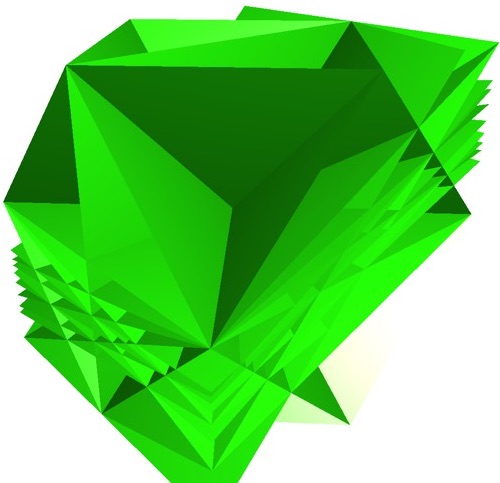}
\end{array}\label{192 picture}
\end{equation}
In fact, there are \emph{two} hyperplane arrangements: one finite arrangement $\scrH$, and the other $\scrH^{\aff}$ is infinite.    When $\scrR$ has only isolated singularities there is a very elegant description of the number of minimal models due to Pinkham \cite{Pinkham}, whereby $\scrH$ is the movable cone.  When $\scrR$ has non-isolated singularities the situation is similar, but mildly more complicated.  This is all recalled in \S\ref{combo section}.

As opined by Reid \cite{Pagoda}, the geometric distinction between isolated and non-isolated cDV singularities is completely artificial.  I am pleased to report that this is as true homologically and noncommutatively as it is geometrically.  With the correct language and concepts, there is little to no distinction between the two cases: all cDV singularities should and can be treated as one.

\subsection{Noncommutative Resolutions and Variants}\label{NC intro sect}
One feature of any cDV singularity $\scrR$ is that all of its birational geometry, which sits `above' $\Spec\scrR$, can be recovered from the category of maximal Cohen--Macaulay modules $\CM\scrR$, in various different ways.  It is worth pausing to explain why this is at least plausible.  Given any minimal model, or more generally any crepant partial resolution $\scrX\to\Spec\scrR$, Van den Bergh \cite{VdB1d} constructs a tilting bundle $\scrO\oplus\scrN$, which after setting $N\colonequals \mathrm{H}^0(\scrN)$ induces a derived equivalence
\begin{equation}
\Db(\coh \scrX)\xrightarrow{\sim}\Db(\mod\End_{\scrR}(\scrR\oplus N))\label{VdBtilt}
\end{equation}
The key fact is that the crepancy of the morphism is equivalent to the homological condition $\End_{\scrR}(\scrR\oplus N)\in\CM\scrR$.  Thus crepancy, a property of the birational geometry `above' $\Spec\scrR$, can be checked on the base ring $\scrR$.

Running the logic backwards, optimistically all $N$ for which $\End_{\scrR}(\scrR\oplus N)\in\CM\scrR$ should be `shadows' of birational geometry.  Remarkably, for cDV singularities this turns out to be true in a very precise way, detailed in \S\ref{AM intro sect} below.  To set language, we say that $M\in\CM\scrR$ is \emph{modifying} if $\End_{\scrR}(R\oplus M)\in\CM\scrR$, and we say it is \emph{maximal modifying} if it is maximal with respect to this property; see \S\ref{NCCR section}.  In that situation, the endomorphism ring $\End_{\scrR}(R\oplus M)\in\CM\scrR$ is called a maximal modification algebra (MMA).

On the other hand a noncommutative crepant resolution (NCCR) is a smooth version of this, which for the sake of this introduction we will take to mean $N\in\CM\scrR$ which is modifying, for which $\End_{\scrR}(\scrR\oplus N)$ has finite global dimension.  In this cDV setting, just as if one minimal model is smooth they all are, if one MMA has finite global dimension, they all do. The terse, and rough, summary of \S\ref{AM section} is the following:
\[
\begin{tikzpicture}[yscale=0.75]
\node (A) at (-1,0) {crepant resolutions};
\node (B) at (-1,-1) {minimal models};
\node (C) at (-1.4,-2) {crepant partial resolutions};
\node (a) at (3,0) {NCCRs};
\node (b) at (3,-1) {MMAs};
\node (c) at (3.9,-2) {modifying algebras};

\draw[<->,densely dotted] (1,0) -- node [above] {$\scriptstyle $} (2,0);
\draw[<->,densely dotted] (1,-1) -- node [above] {$\scriptstyle $} (2,-1);
\draw[<->,densely dotted] (1,-2) -- node [above] {$\scriptstyle $} (2,-2);
\end{tikzpicture}
\]
where if the top line exists, then the first line is equal to the second line.  The point of course is that NCCRs (like crepant resolutions) do not in general exist, and so the objects we aim for are MMAs (like minimal models).  We view the case when minimal models are smooth, respectively MMAs have finite global dimension, as nothing more than a happy coincidence.  In this way, all cDV singularities get treated at once, avoiding arbitrary dividing lines.  In fact, restricting to minimal models in itself a touch artificial: most, if not all, results for minimal models below have counterparts for the more general crepant partial resolutions.
  
The dotted lines above are still a bit vague, since in order to obtain more precise relationships, and bijections, we need to specify how we count, and restrict the generality on the noncommutative side.

\subsection{Viewpoint from Cluster Theory}
At the level of functors and derived categories, it turns out that Bridgeland's flop functor \cite{Bridgeland} is, suitably interpreted, a `cluster' mutation \cite{IW3, HomMMP}.  Thus the above birational geometry can be achieved using various concepts that arise in the (additive) categorification of cluster algebras, which we now outline.  One of the reasons this is beneficial is that cluster theory is designed for iterating mutations, whereas iterating flops in birational geometry is usually more tricky. Another is that it gives rise to the idea that the cluster tilted algebra (aka the contraction algebra, below) should be the classifying object.

Consider, as in \S\ref{NC intro sect}, the category of maximal Cohen--Macaulay modules $\CM\scrR$ where $\scrR$ is cDV.  When $\scrR$ is isolated the \emph{stable} category $\uCM\scrR$ is both Hom-finite and 2-CY, two of the key properties shared with categories that arise in cluster theory.  Alas, when $\scrR$ is not isolated, both properties fail, but they can be recovered after some suitable reinterpretation.  In either case, being Krull--Schmidt and appropriately 2-CY, the category $\uCM \scrR$ serves as a `model' of various versions of cluster theory.

Again, there is a plethora of related concepts, rigorously introduced in \S\ref{NCCR section}, but for now are roughly summarised by the following diagram
\[
\begin{tikzpicture}[yscale=0.75]
\node at (-1,1) {$\scrR$ isolated};
\node at (3,1) {$\scrR$ general};
\draw[densely dotted] (-3,0.5)--(6.5,0.5);
\node (A) at (-1,0) {cluster tilting objects};
\node (B) at (-1,-1) {maximal rigid objects};
\node (C) at (-1,-2) {rigid objects};
\node (a) at (3.1,0) {CT objects};
\node (b) at (4.3,-1) {maximal modifying objects};
\node (c) at (3.6,-2) {modifying objects};

\draw[->,densely dotted] (1,0) -- node [above] {$\scriptstyle $} (1.8,0);
\draw[->,densely dotted] (1,-1) -- node [above] {$\scriptstyle $} (1.8,-1);
\draw[->,densely dotted] (1,-2) -- node [above] {$\scriptstyle $} (1.8,-2);
\end{tikzpicture}
\]
Again, there are two main points: (1) the left/right dichotomy is artificial, we should always work generally and so in the right hand column, and (2) the top line need not exist, but if it does, the top line equals the middle line.  Regardless, we \emph{aim} for the middle line, and view the occasions when the top line exists as a happy coincidence.

The analogy between birational geometry and cluster theory turns out to be very useful, and this will be expanded upon in later sections.  Analogies, however, only take you so far, and there are several key distinguishing features of this geometric cDV setting.  The easiest, namely that cluster tilting objects need not exist, is solved by instead studying maximal rigid (actually, maximal modifying) objects.  Algebraically, that is now a fairly standard viewpoint.  However, more fundamentally, the quivers of our endomorphism algebras \emph{usually} have loops, and \emph{always} have 2-cycles, a fact which still seems to stir up remarkably strong emotions.  Really, loops and 2-cycles are not a big deal.  Their existence just means that the theory is not \emph{entirely} combinatorial.  It would be much more perverse if a huge chunk of the birational geometry of $3$-folds was controlled by extremely elementary combinatorial rules.  Birational geometry is just more complicated than that.

\subsection{Auslander--McKay Correspondence}\label{AM intro sect}
For Kleinian singularities, the McKay correspondence \cite{McKay}, as reformulated by Auslander \cite{Aus2,Aus3}, gives a bijection between indecomposable non-free Cohen--Macaulay modules and exceptional curves in the minimal resolution. It turns out that this can be lifted to cDV singularities, at the cost of introducing \emph{two} bijections, needed to incorporate the added complexities of $3$-fold geometry.

First, there is a bijection between maximal modifying (=MM) $\scrR$-module generators and minimal models, and then for each such pair there is a further bijection (in parts \eqref{NCvsC min model intro 1} and \eqref{NCvsC min model intro 2} below) which is very similar to the classical Auslander--McKay Correspondence for surfaces.  Part \eqref{NCvsC min model intro 3} describes how these two bijections relate.

\begin{thm}[Auslander--McKay Correspondence in Dimension $3$]\label{AM intro}
There is a bijection
\[
\begin{array}{c}
\begin{tikzpicture}
\node (A) at (-1,0) {$\{ \mbox{basic \textnormal{MM} $\scrR$-module generators}\}$};
\node (B) at (6,0) {$
\{\mbox{minimal models $f_i\colon \scrX_i\to\Spec \scrR$} \}
$};
\draw[<->] (1.85,0) -- node [above] {$\scriptstyle $} (3,0);
\node at (0,0.175) {$\phantom -$};
\end{tikzpicture}
\end{array}
\]
where the left-hand side is taken up to isomorphism, and the right-hand side is taken up to isomorphism of the $\scrX_i$ compatible with the morphisms $f_i$.  Under this bijection:
\begin{enumerate}
\item\label{NCvsC min model intro 1} For any \textnormal{MM} generator, its non-free indecomposable summands are in one-to-one correspondence with the exceptional curves in the corresponding minimal model.
\item\label{NCvsC min model intro 2} For any fixed  \textnormal{MM} generator $N$, the quiver of $\uEnd_\scrR(N)$ encodes the dual graph of the corresponding minimal model.
\item\label{NCvsC min model intro 3} The mutation graph of the  \textnormal{MM} generators coincides with the flops graph of the minimal models.
\end{enumerate}
\end{thm}

All undefined terminology, and detailed descriptions of the bijections, are given in \S\ref{NCCR section}. Most of the above does not require the restriction to minimal models, and indeed this is needed for the applications to autoequivalences and more general terminal singularities.

The graphs in \eqref{NCvsC min model intro 3} are simply the framework to express the relationship between flops and mutation on a \emph{combinatorial} level.  The real power comes from lifting this to the level of functors, as this then impacts on autoequivalences, group actions and stability conditions.  All these rely on an understanding of the tilting theory of the associated MMAs, which we sketch in \S\ref{tilt sect intro}.  The second piece of extra homological information comes from the left hand side of the bijection: relaxing that to include \emph{all} MM modules (not just the generators) gives rise to the \emph{affine action} on the derived category sketched in \S\ref{auto intro}, and unveils t-structures that the birational geometry does not `see'.

\subsection{Contraction Algebras}
Before proceeding to categorical and functorial considerations, we briefly outline the role of the conjectural classifying object, namely the contraction algebra, and how it is constructed.

First, for an arbitrary crepant partial resolution $f\colon\scrX\to\Spec\scrR$ with $\scrR$ cDV, on one hand the Van den Bergh equivalence \eqref{VdBtilt} equips us with $N\in\CM\scrR$, and so we may form the stable endomorphism algebra  $\uEnd_\scrR(N)$.  On the other hand, over the unique closed point $p$ of $\Spec\scrR$, consider $\mathrm{C}\colonequals f^{-1}(p)$ equipped with its reduced scheme structure, namely $\mathrm{C}^{\rm red}$.  The idea is to construct, using noncommutative deformation theory,
\[
\begin{tikzpicture}
\node (A) at (0,0) {$\mathrm{C}^{\rm red}$ in $\scrX$ \kern 4pt };
\node (B) at (5.5,0) {\kern 4pt an algebra $\CA$.};
\draw [->,decorate, 
decoration={snake,amplitude=.4mm,segment length=3mm,post length=1mm}]  (A) to node[above] {\scriptsize (pro)representing object} (B);
\end{tikzpicture}
\]
The \emph{contraction  algebra} $\CA$ is an analytic invariant of the point  $p$, but in the algebraic setting it is satisfyingly possible to glue all the above local constructions together, and so obtain a sheaf of algebras $\scrD$, such that $\scrD_p$ is morita equivalent to $\CA$, for each point $p$.  Locally, it turns out that $\CA\cong \uEnd_\scrR(N)$, and so the algebra constructed using deformation theory is the same as that from cluster theory.  This has many consequences.

The details of the construction, using NC deformation theory, are given in \S\ref{def theory section}.  Perhaps the key point is that commutative deformation theory can be recovered from this, by losing information, which is why $\CA$ recovers all known curve invariants, whilst containing \emph{extra} information in the form of an algebra structure.  As a rough rule of thumb, noncommutative deformations of the fixed sheaf $\scrO_{\mathrm{C}^{\rm red}}$ detects how the curve $\mathrm{C}^{\rm red}$ moves, \emph{and} all of its multiples (see \S\ref{sec: curves count} below).  Thus, the contraction algebra bundles all this numerical information into one local object, together with an algebra structure.

The basic philosophy is that the contraction algebra controls all the birational geometry locally around $p$, and classifies in the situation when $\scrX$ is smooth.  Indeed, the Contraction Theorem, stated in \ref{contraction theorem}, asserts that the difference between the cases \ding{192} and \ding{193} in \eqref{two cases intro} can be detected by the finite dimensionality (or otherwise) of the contraction algebra.  The Classification Conjecture, on the other hand, asserts that the isomorphism class of the contraction algebra of a smooth minimal model analytically locally classifies the cDV singularity.  The precise statement (see \ref{class:conj}) is mildly more involved, given minimal models are not unique, but it can still be reduced to an isomorphism problem.

It is worth pausing to unpack the above paragraph.  A particular special case is that smooth, irreducible, $3$-fold flopping contractions are conjecturally classified by finite dimensional data, in the form of $\CA$ and its algebra structure.  Via various translations through superpotential algebras (see e.g.\ \cite[3.1]{BW}), this means that smooth irreducible flops are conjecturally classified by certain elements in the noncommutative formal power series ring $\mathbb{C}\llangle x,y\rrangle$.

\subsection{Curve Counting}\label{sec: curves count}
Much of the power of contraction algebras comes from the fact that they admit three distinct but compatible structures, explained in \S\ref{main contalg subsection}.  In short, $\CA$ can be viewed as a sheaf of $\scrX$ (via \eqref{VdBtilt}), as a module on $\scrR$, or as an algebra.  

When the support of $\CA$, as an $\scrR$-module, is finite, in particular $\CA$ is finite dimensional as a vector space.  Through this, we can thus associate a new finite-dimensional invariant to every flopping (and more generally flipping!) contraction.  It is possible to then lose information in various ways, by taking the dimension of  $\CA$, or by considering its abelianisation, or its quiver, and from this recover known invariants.  One such example of this is Toda's dimension formula \cite{TodaGV}, which in the case of an smooth irreducible flopping contraction states that
\[
\dim_{\mathbb{C}}\CA = \underbrace{\dim_{\mathbb{C}} \CA^{\rm ab}}_{n_1} + \sum_{i=2}^{\ell} n_i\cdot i^2,
\]
where $\ell$ is the \emph{length}, which is determined by the elephant, and \emph{all} $n_i> 0$.  The $n_i$ are called the Gopakumar--Vafa (GV) invariants, and are defined via deforming the flopping curve into Atiyah flops (see \S\ref{Toda section}).    We remark that the GV invariants are a property of the isomorphism class of $\CA$, but it is still an open problem to extract them easily.

\subsection{Tilting}\label{tilt sect intro}
Whilst the stable endomorphism ring $\uEnd_{\scrR}(N)$, in the guise of the contraction algebra, controls and conjecturally classifies the local geometry, the technical backbone of the homological algebra is obtained through studying the more standard endomorphism ring $\End_\scrR(N)$, and in particular its tilting theory.

To ease notation set $\Lambda=\End_{\scrR}(N)$. Recall that by a classical tilting module we mean $T\in\mod\Lambda$ such that (1) $\pd_{\Lambda}T\leq 1$, (2) $\Ext^1_{\Lambda}(T,T)=0$, and (3) there exists an exact sequence $0\to\Lambda\to T_0\to T_1\to 0$ with each $T_i\in\add T$.  This last condition is really just a generation condition.   Our motivation for studying such objects arises from the fact each such $T$ induces a derived equivalence
\[
\Db(\mod\End_{\scrR}(N))\to\Db(\mod\End_{\Lambda}(T))
\]
and furthermore $\End_{\Lambda}(T)\cong\End_{\scrR}(L)$ for some modifying $\scrR$-module $L$.  Thus, first controlling and varying these $T$, then composing the resulting equivalences gives derived autoequivalences of any fixed $\End_{\scrR}(N)$. Hence, via  \eqref{VdBtilt}, this generates autoequivalences of $\Db(\coh\scrX)$.

Recall from \S\ref{geo overview intro} that associated to $\scrX\to\Spec\scrR$ is a finite hyperplane arrangement $\scrH$ and an infinite version $\scrH^{\aff}$.  Perhaps the key technical result, which we only sketch here, is that there are a series of maps
\[
\begin{tikzpicture}
\node (A) at (0,0) {$\mathrm{MM}\,\scrR$};
\node (B) at (4,0) {$\mathrm{reftilt}\End_{\scrR}(N)$};
\node (C) at (8,0) {$\mathrm{tilt}\,e_{\scrJ}\Uppi e_{\scrJ}$};
\node (D) at (11,0) {$\mathsf{Cham}\,\scrH^{\aff}$};
\draw[->] (A) -- node[above] {$\scriptstyle \Hom_{\scrR}(N,-)$} node[below] {$\scriptstyle \sim$}(B);
\draw[right hook->] (B) -- node[above] {$\scriptstyle (\scrR/g)\otimes_{\scrR}-$}  (C);
\draw[->] (C) -- node[above] {$\scriptstyle \sim$} (D);
\end{tikzpicture}
\]
where the outers are always bijections, and the middle is an injection.  When $\scrR$ is isolated, all are bijections.  In the above notation, $\mathrm{MM}\,\scrR$ is the set of maximal modifying $\scrR$-modules, $\mathrm{reftilt}\End_{\scrR}(N)$ is a certain subset of tilting modules, and $e_{\scrJ}\Uppi e_{\scrJ}$ is a contracted preprojective algebra, which is `surfaces data'.  Thus the above maps link commutative $3$-fold objects with noncommutative $3$-dimensional algebras, with noncommutative $2$-dimensional algebras, and lastly with the combinatorics of $\scrH^{\aff}$.  

There is a similar series of bijections for the finite arrangement $\scrH$, if we restrict each set appropriately, and there are versions of the above bijections for arbitrary (non-maximal) modifying modules too.  This generality is needed for applications to flopping contractions $\scrX\to\Spec\scrR$ where $\scrX$ is mildly singular, but can otherwise be ignored.

Ultimately, there are two key points. First, there is an operation of mutation on the first three sets, and wall-crossing on the last, and the maps intertwine all these actions.  This means that usually complicated mutation behaviour is modelled on unusally simple wall crossing combinatorics.  Second, the main point is that the maps (which are all bijections when $\scrR$ is isolated) exist for the infinite arrangement $\scrH^{\aff}$, and not just the finite movable cone $\scrH$.  Hence, we obtain not just the `birational geometry' action by the flop functors below, we obtain a much larger action, one that does not yet have a birational geometry interpretation.

\subsection{Autoequivalences and Stability Conditions}\label{auto intro}
Given any flopping contraction $\scrY\to\scrY_{\con}$ between quasi-projective (e.g.\ projective!) $3$-folds, where $\scrY$ has at worst isolated cDV singularities, then the above complete-local contraction algebra technology has global consequences.

Firstly, there exists a functorial triangle
\begin{equation}
\RHom_{\scrY}(\scrE,-)\otimes_{\CA}^{\bf L}\scrE\to \Id\to\twistGen_{}\to\label{sph disp intro}
\end{equation}
where $\scrE$ is the universal bundle of the NC deformation theory, viewed as a sheaf on $\scrY$.  The third term $\twistGen_{}$ should be interpreted as a twist over a spherical functor, with base $\Db(\mod\CA)$.  Variations of the above are possible, and indeed desirable: choosing only a subset $\scrJ$ of the flopping curves (which need not contract algebraically), we obtain a different twist functor $\twistGen_{\scrJ}$, this time with base given by noncommutative deformations of the chosen subset of the curves.  More details are given in \S\ref{autostab section}.

We remark here that autoequivalences are one situation where the difference between isolated and non-isolated cDV singularities becomes stark.  Flopping contractions, to isolated singular points, always give rise to global autoequivalences.  In contrast, divisor-to-curve contractions sometimes do, and sometimes do not \cite[\S 6]{DW4}.  On the positive side, at least this latter fact is consistent across mirror symmetry (see \S\ref{HMS section}).

\medskip
Returning to the global flopping contraction $\scrY\to\scrY_{\con}$, the twist autoequivalences described above \emph{don't} braid, they satisfy the relations of the \emph{pure} braid group.  The presentation of such groups is ghastly, and so we are forced to rely on their description as fundamental groups of complexified complements $\scrH$ and $\scrH^{\aff}$ from \S\ref{geo overview intro}; this is recapped in \S\ref{autostab section}.  Writing $\scrZ$ for the complexified complement of $\scrH$, and  $\scrZ_{\aff}$ for the $\scrH^{\aff}$ version, the upshot is that the above autoequivalences knit together, at least for local flops $\scrX\to\Spec\scrR$, to  form the following group homomorphisms
\[
\begin{tikzpicture}
\node (A1) at (0,0) {$\uppi_1(\scrZ)$};
\node (A2) at (0,-1.5) {$\uppi_1(\scrZ_{\aff})$};
\node (B) at (2.5,0) {$\mathrm{Auteq}\Db(\coh \scrX)$};
\draw[->] (A1) -- node[above] {$\scriptstyle\upvarphi$}(B);
\draw[->] (A2) --  node[gap]  {$\scriptstyle\widetilde{\upvarphi}$}(B);
\draw[->] (A1) -- (A2);
\end{tikzpicture}
\]
There are versions of the above in the global $\scrY\to\scrY_{\con}$ case, and also versions for arbitrary cDV singularities. However, for non-isolated cDV singularities there are repeats in the chambers, and this repetition creates difficulties in stating things elegantly. The same problem arises for derived categories of contraction algebras even in the isolated case, and this is explained in detail in \cite[7.3]{August2}.

In the easiest case, when $\scrH$ is a Type $A$ root system, $\uppi_1(\scrZ)$ is the classical pure braid group, and $\uppi_1(\scrZ_{\aff})$ is its affine version.  Hence we lazily refer to the above as the `pure braid action' and the `affine pure braid action', regardless of whether $\scrH$ is Coxeter.  

It is known that the homomorphism $\upvarphi$ is injective, whilst $\widetilde{\upvarphi}$ is only conjectured to be so.  The image of $\upvarphi$ is important from a birational geometry viewpoint, being the group generated by all flop functors and their inverses.  On the other hand, the images of both  $\upvarphi$ and $\widetilde{\upvarphi}$ play a key role in the description of Bridgeland stability conditions in this setting, namely as the Galois groups of regular covering maps explained in \S\ref{autostab section}.

\subsection{Viewpoint from Mirror Symmetry}\label{HMS section}
Given a crepant resolution $\scrX\to\Spec\scrR$, where $\scrR$ is cDV, it is of value to both sides of mirror symmetry to understand and construct the symplectic mirror $X$ to the smooth CY $3$-fold $\scrX$.  Aside from some embarrassingly easy examples of $\scrR$, which includes the Atiyah flop, some toric examples, and not much else besides, this problem is wide open.

It is instructive to take one step back, and explain first the key behaviour that is required to be mirrored. In the case of the Atiyah flop, which is the only example of a flopping $(-1,-1)$-curve in birational geometry, the structure sheaf $\mathsf{E}=\scrO_{\mathrm{C}^{\rm red}}$ is a spherical object, namely
\[
\Ext^i_{\scrX}(\mathsf{E},\mathsf{E})=
\begin{cases}
\mathbb{C}&\mbox{if }i=0,3\\
0&\mbox{else.}
\end{cases}
\] 
The functorial triangle  \eqref{sph disp intro} in this case dramatically reduces in complexity: here $\CA=\mathbb{C}$, $\scrE=\mathsf{E}$, and so the autoequivalence obtained is just the standard spherical twist, which should be viewed as the analogue of the symplectic Dehn twist of Seidel \cite{ST}.

\medskip
It is somewhat inconvenient that the Atiyah flop is the \emph{only} example for which $\mathsf{E}=\scrO_{\mathrm{C}^{\rm red}}$ is spherical.  For all other single-curve flops, either
\begin{equation}
\Ext^i_{\scrX}(\mathsf{E},\mathsf{E})
=
\begin{cases}
\mathbb{C}&\mbox{if }i=0,3\\
\mathbb{C}&\mbox{if }i=1,2\\
0&\mbox{else}
\end{cases}
\quad
\mbox{or}
\quad
=
\begin{cases}
\mathbb{C}&\mbox{if }i=0,3\\
\mathbb{C}^2&\mbox{if }i=1,2\\
0&\mbox{else.}
\end{cases}\label{single curve cohomology}
\end{equation}
Again, don't be fooled: there is only one family of single-curve flops, the so-called Pagoda family, which satisfy the left hand conditions in \eqref{single curve cohomology}.  All other flops are hidden under the umbrella of the right hand conditions.  These are the flopping curves with normal bundle $(-3,1)$.

What were previously spheres now behave more like $S^1\times S^2$, Lens spaces or Seifert fibre spaces.  The upshot of \S\ref{auto intro} is that the behaviour of these objects, and the symmetries that they generate, are \emph{not} governed just by their cohomology: the higher $\mathrm{A}_\infty$-operations now play a much more central role. Indeed, the associated autoequivalences are controlled by the noncommutative deformations $\scrE$ of $\mathsf{E}$, which in turn are determined by the $\mathrm{A}_\infty$-algebra $\Ext_{\scrX}^*
(\mathsf{E},\mathsf{E})$; see \S\ref{NC multiple sction} for full details.

\medskip
Thus, at a bare minimum, to mirror a single-curve flop  requires the existence of a Lagrangian $\mathsf{E}$ inside $X$ whose Floer cohomology satisfies one of the cases in \eqref{single curve cohomology}.  The difficultly, exactly as in the algebraic geometric setting, is finding objects $\mathsf{E}$ that not only satisfy the right hand conditions in \eqref{single curve cohomology}, but \emph{furthermore} whose $\mathrm{A}_\infty$-operations
\[
\mathsf{m}_i\colon \Ext^1_{\scrX}(\mathsf{E},\mathsf{E})^{\otimes i}\to \Ext^2_\scrX(\mathsf{E},\mathsf{E})
\]
after dualizing furthermore present an algebra
\[
\CA\cong \frac{\mathbb{C}\langle\!\langle\, \Ext^1_\scrX(\mathsf{E},\mathsf{E})^*\rangle\!\rangle}{\mathsf{m}^\vee(\Ext^2_\scrX(\mathsf{E},\mathsf{E})^*)}
\]
which is \emph{finite dimensional} as a vector space.  In the case of the right hand conditions in \eqref{single curve cohomology}, this is asking us to realise various finite dimensional quotients of $\mathbb{C}\langle\!\langle x,y\rangle\!\rangle$ as examples in both algebraic and symplectic geometry.  The case of the left hand conditions in \eqref{single curve cohomology} is described in \cite{MakWu}.

\medskip
If we remove the requirement of a flopping contraction from the above discussion, and instead ask to mirror an arbitrary crepant resolution $\scrX\to\Spec\scrR$ with $\scrR$ cDV, things become more complicated.  One property exhibited even by $S^1\times S^2$-type objects, which precisely mirrors one of the problems described in \S\ref{auto intro}, is that it is sometimes not possible to extend symmetries (or autoequivalences) associated to such objects globally.  Indeed, if $L = S^1\times S^2 \subset X$ is a Lagrangian, then there is a non-compactly-supported symplectomorphism  $\mathrm{id} \times \uptau_{S^2} \in \mathrm{Symp}(T^*S^1\times T^*S^2)$ which, alas, need not extend to a global symplectomorphism of $X$ \cite[4.12]{SW}.  Sometimes it does, sometimes it does not.

\subsection{Algebraic vs Complete Local} 
Much of the frustration levelled at the theory of cDV singularities stems from results being stated in the complete local setting, when ultimately most people want to work algebraically.  Kleinian singularities are global!  Give us good old algebra!  Here, we really do have to be careful what we wish for.  

The main issue is that cDV singularities do not admit a grading, in general.  They are just fundamentally \emph{more complicated} than Kleinian singularities. By all means, working Zariski locally is possible, and most statements in this survey have an algebraic counterpart (see e.g.\ \cite[\S6]{DW1}, or \cite[\S2]{IW3}).  But you have to question whether you want to be the person who counts things up to additive equivalence, defines mutation up to Morita equivalence, and who revels in the use of lifting numbers to skirt round issues involving the failure of Krull--Schmidt.

Consequently, this survey restricts exclusively to the complete local setting. The technical advantages this affords means that results can be stated much more simply.  The Zariski-local statements often mask otherwise simple concepts, and have a cumulative draining effect: the curious and the brave are invited to read \cite[\S2]{IW3} to witness first hand this draining phenomenon.

\subsection*{Conventions} 
We work over the field $\mathbb{C}$ throughout, and at various points this is very important.  When composing arrows, functors or morphisms, any product not containing the $\circ$ symbol, such as $ab$, will mean $a$ followed by $b$.

\subsection*{Acknowledgements}  
I place on record my thanks to all my collaborators and PhD students who have contributed to the mathematics from which this survey has sprung. Particularly  Osamu Iyama for entertaining the idea, over 10 years ago now, that birational geometry is controlled by generalisations of maximal rigid objects, and Will Donovan for relentlessly, but patiently, explaining that mathematical communication might actually be improved if we purge all passively aggressive remarks from our papers.  Who knew!  Thanks are also due to Jenny August, Gavin Brown, Yuki Hirano and Ivan Smith for collaborations and comments, and to Ali Craw, Okke van Garderen, Ben Davison, Miles Reid for various and wider discussions over the years.

\subsection*{Other Surveys}  
A more algebraic introduction to some of the ideas presented here can be found in the PhD thesis of August \cite[\S1--3]{AugustThesis}, which focuses on the derived category, cluster theory and contraction algebra viewpoints for isolated cDV singularities.  More general surveys on noncommutative resolutions include \cite{LeuschkeSurvey,WemyssSurvey}, but both are somewhat out of date. The best geometric introduction to cDV singularities remains the original work of Reid \cite{Pagoda,YPG,ReidCanonical}, to which much intellectual debt is owed.

\section{Birational Geometry}\label{geometry section}

It is impossible to convey the role of cDV singularities without the language of the MMP.  The reader is referred to \cite{Pagoda, KollarMori} for a general overview; we maintain a specific focus on the aspects most relevant to cDV singularities here.  As in the introduction, let $\scrR$ denote a complete local cDV singularity.   Automatically $\scrR$ is a normal domain.

\subsection{Generalities and Minimal Models}\label{min models sect} Recall that a normal scheme $\scrX$ is said to be $\mathds{Q}$-factorial if for every Weil divisor $D$, there exists $n\in\mathbb{N}$ for which $nD$ is Cartier. Further, a projective birational morphism $f\colon \scrX\to \Spec\scrR$ is called crepant if $f^*\omega_{\scrR}=\omega_\scrX$. A $\mathds{Q}$-factorial terminalisation, or minimal model, of $\Spec\scrR$ is a crepant projective birational morphism $f\colon \scrX\to \Spec\scrR$ such that $\scrX$ has only $\mathds{Q}$-factorial terminal singularities (see \S\ref{terminal}). When $\scrX$ is furthermore smooth, we call $f$ a crepant resolution.  

\begin{example}\label{toric ex}
We give two easy, toric, examples.
\begin{enumerate}
\item\label{toric ex 1}  The conifold, or Atiyah flop, given by $\scrR=\mathbb{C}\llsq u,v,x,y\rrsq/(uv-xy)$.  Blowing up the ideal $(u,x)$ gives rise to $\scrX\to\Spec\scrR$, and blowing up $(u,y)$ gives rise to $\scrX^+\to\Spec\scrR$.  Both $\scrX$ and $\scrX^+$ are smooth, and the normal bundle of the exceptional curves in both cases is $\scrO_{\mathbb{P}^1}(-1)^{\oplus 2}$.  

\item\label{toric ex 2} The suspended pinch point, given by $\scrR=\mathbb{C}\llsq u,v,x,y\rrsq/(uv-x^2y)$.  This is singular along the $x$-axis.  There are precisely three crepant resolutions, sketched below.

\[
\begin{array}{c}
\begin{tikzpicture}
\node at (-3,0) {\begin{tikzpicture}[scale=0.6]
\coordinate (T1) at (1.7,2.5);
\coordinate (T) at (1.9,1.8);
\coordinate (TM) at (1.92,1.7);
\coordinate (B) at (2.1,0.9);
\foreach \y in {0,0.1,0.2,...,1}{ 
\draw[line width=0.5pt,red!50] ($(T)+(\y,0)$) to [bend left=25] ($(B)+(\y,0)$);
\draw[line width=0.5pt,red!50] ($(T)+(-\y,0)$) to [bend left=25] ($(B)+(-\y,0)$);;}
\draw[red,line width=1pt] (T1) to [bend left=25] (TM);
\draw[line width=1pt,red] (T) to [bend left=25] (B);
\draw[color=blue!60!black,rounded corners=6pt,line width=0.75pt] (0.5,0) -- (1.5,0.3)-- (3.6,0) -- (4.3,1.5)-- (4,3.2) -- (2.5,2.7) -- (0.2,3) -- (-0.2,2)-- cycle;
\end{tikzpicture}};
\node at (0,0) {\begin{tikzpicture}[scale=0.6]
\coordinate (T) at (1.9,2);
\coordinate (TM) at (2.12-0.05,1.5-0.1);
\coordinate (BM) at (2.12-0.09,1.5+0.12);
\coordinate (B) at (2.1,1);
\foreach \y in {0.1,0.2,...,1}{ 
\draw[line width=0.5pt,red!50] ($(T)+(\y,0)+(0.02,0)$) to [bend left=25] ($(B)+(\y,0)+(0.02,0)$);
\draw[line width=0.5pt,red!50] ($(T)+(-\y,0)+(-0.02,0)$) to [bend left=25] ($(B)+(-\y,0)+(-0.02,0)$);;}
\draw[red,line width=1pt] (T) to [bend left=25] (TM);
\draw[red,line width=1pt] (BM) to [bend left=25] (B);
\draw[color=blue!60!black,rounded corners=6pt,line width=0.75pt] (0.5,0) -- (1.5,0.3)-- (3.6,0) -- (4.3,1.5)-- (4,3.2) -- (2.5,2.7) -- (0.2,3) -- (-0.2,2)-- cycle;
\end{tikzpicture}};
\node at (3,0) {\begin{tikzpicture}[scale=0.6]
\coordinate (T) at (1.9,2.2);
\coordinate (B) at (2.1,1.3);
\coordinate (BM) at (2.07,1.45);
\coordinate (B1) at (2.3,0.65);
\foreach \y in {0,0.1,0.2,...,1}{ 
\draw[line width=0.5pt,red!50] ($(T)+(\y,0)$) to [bend left=25] ($(B)+(\y,0)$);
\draw[line width=0.5pt,red!50] ($(T)+(-\y,0)$) to [bend left=25] ($(B)+(-\y,0)$);;}
\draw[red,line width=1pt] (BM) to [bend left=25] (B1);
\draw[line width=1pt,red] (T) to [bend left=25] (B);
\draw[color=blue!60!black,rounded corners=6pt,line width=0.75pt] (0.5,0) -- (1.5,0.3)-- (3.6,0) -- (4.3,1.5)-- (4,3.2) -- (2.5,2.7) -- (0.2,3) -- (-0.2,2)-- cycle;
\end{tikzpicture}};
\node at (0,-2) {\begin{tikzpicture}[scale=0.6]
\draw [red] (1.1,0.75) -- (3.1,0.75);
\draw[color=blue!60!black,rounded corners=6pt,line width=0.75pt] (0.5,0) -- (1.5,0.15)-- (3.6,0) -- (4.3,0.75)-- (4,1.6) -- (2.5,1.35) -- (0.2,1.5) -- (-0.2,1)-- cycle;
\end{tikzpicture}};
\draw[->, color=blue!60!black] (0,-1) -- (0,-1.5);
\draw[->, color=blue!60!black] (-2,-1) -- (-1.5,-1.5);
\draw[->, color=blue!60!black] (2,-1) -- (1.5,-1.5);
\end{tikzpicture}
\end{array}
\]
\end{enumerate}
\end{example}

\begin{example}\label{nontoric ex}
We give three less easy examples.  More are given in \ref{cut ex}\eqref{cut ex 2} and \S\ref{Type A section}.
\begin{enumerate}
\item\label{nontoric ex 1} The $cA_2$ singularity given by  $\scrR=\mathbb{C}\llsq u,v,x,y\rrsq/(uv-(x^2+y^3)x)$. Blowing up $(u,x)$ gives $\scrX\to\Spec\scrR$, and blowing up $(u,x^2+y^3)$ gives $\scrX^+\to\Spec\scrR$.  Neither $\scrX$ or $\scrX^+$ is smooth, they both have a unique singular point. Pictorially:
\[
\begin{tikzpicture}
\node at (-1.5,0) {$\begin{tikzpicture}[scale=0.5]
\coordinate (T) at (1.9,2);
\coordinate (B) at (2.1,1);
\draw[red,line width=1pt] (T) to [bend left=25] (B);
\filldraw[blue!60!black] (B) circle (2pt); 
\draw[color=blue!60!black,rounded corners=5pt,line width=1pt] (0.5,0) -- (1.5,0.3)-- (3.6,0) -- (4.3,1.5)-- (4,3.2) -- (2.5,2.7) -- (0.2,3) -- (-0.2,2)-- cycle;
\end{tikzpicture}$};
\node at (1.5,0) {$\begin{tikzpicture}[scale=0.5]
\coordinate (T) at (2.1,2);
\coordinate (B) at (1.9,1);
\draw[red,line width=1pt] (T) to [bend right=25] (B);
\filldraw[blue!60!black] (T) circle (2pt); 
\draw[color=blue!60!black,rounded corners=5pt,line width=1pt] (0.5,0) -- (1.5,0.3)-- (3.6,0) -- (4.3,1.5)-- (4,3.2) -- (2.5,2.7) -- (0.2,3) -- (-0.2,2)-- cycle;
\end{tikzpicture}$};
\node at (0,-2) {\begin{tikzpicture}[scale=0.5]
\filldraw [red] (2.1,0.75) circle (1pt);
\draw[color=blue!60!black,rounded corners=5pt,line width=1pt] (0.5,0) -- (1.5,0.15)-- (3.6,0) -- (4.3,0.75)-- (4,1.6) -- (2.5,1.35) -- (0.2,1.5) -- (-0.2,1)-- cycle;
\end{tikzpicture}};
\draw[->, color=blue!60!black] (-0.75,-1) -- (-0.25,-1.5);
\draw[->, color=blue!60!black] (0.75,-1) -- (0.25,-1.5);

\end{tikzpicture}
\]
where the dots on the exceptional curves are, complete locally, the singular point $uv=x^2+y^3$, which is factorial.  Hence $\scrX$ and $\scrX^+$ are minimal models of $\Spec\scrR$.

\item  The $cA_3$ singularity given by  $\scrR_\uplambda=\mathbb{C}\llsq u,v,x,y\rrsq/(uv-xy(x+y)(x+\uplambda y))$, where $\uplambda\in\mathbb{C}$ with  $\uplambda\neq 0,1$.  There are infinitely many non-isomorphic $\scrR_\uplambda$, which demonstrates that cDV singularities have moduli.  Regardless, each $\scrR_\uplambda$ has 24 crepant resolutions, described in \S\ref{min model A sect}.  

\item\label{nontoric ex 2} \cite[2.25]{DW4} The singularity $\scrR = \mathbb{C}\llsq x, y, z, t\rrsq/(x^3 - xyt - y^3 + z^2)$, which is singular along the $t$-axis.  This is $cD_4$ at the origin.  Blowing up the singular locus $(x,y,z)$ gives rise to a crepant resolution, which pictorially is \ding{193} in \eqref{two cases intro}.
\end{enumerate}
\end{example}

We next describe flops and flopping curves, but this first requires a brief discussion on contracting curves generally. One of the benefits of the complete local setting, over the algebraic setting, is that we can contract curves at will.  If $f\colon \scrX\to \Spec \scrR$ is a projective birational morphism with $\RDerived f_*\scrO_\scrX=\scrO_\scrR$ and at most one-dimensional fibres, write $\bigcup_{i=1}^n\mathrm{C}_i$ for the fibre above the origin, with reduced scheme structure, so that each $\mathrm{C}_i\cong\mathbb{P}^1$. For any subset $I$ of  these curves, since $\scrR$ is complete local we may factorise $f$ as
\[
\scrX\xrightarrow{g} \scrX_{\con}\xrightarrow{h} \Spec \scrR
\]
where $g$ contracts $\mathrm{C}_i$ to a closed point if and only if $i\in I$, and further both $\RDerived g_*\scrO_\scrX=\scrO_{\scrX_{\con}}$ and $\RDerived h_*\scrO_{\scrX_{\con}}=\scrO_{\scrR}$.

In the contracting morphism $g$, more than just $\bigcup_{i\in I}\mathrm{C}_i$ in $\scrX$ may be contracted.  The situation when \emph{only} this is contracted, namely $g$ is an isomorphism away from $\bigcup_{i\in I}\mathrm{C}_i$, is of particular interest.  To set language, recall that a $\mathds{Q}$-Cartier divisor $D$ is called $g$-nef if $D\cdot \mathrm{C}\geq 0$ for all curves contracted by $g$, and $D$ is called $g$-ample if $D\cdot \mathrm{C}> 0$ for all curves contracted by $g$. 

\begin{defin}\label{flopsdefin}
With $f\colon \scrX\to \Spec \scrR$ as above, choose $\bigcup_{i\in I}\mathrm{C}_i$ in $\scrX$, contract them to give $g\colon \scrX\to \scrX_{\con}$, and suppose that $g$ is an isomorphism away from $\bigcup_{i\in I}\mathrm{C}_i$.  Then we say that $g^+\colon \scrX^+\to \scrX_{\con}$ is the flop of $g$ if for every line bundle $\scrO(D)$ on $\scrX$ such that $-D$ is $g$-nef, then the proper transform of $D$ is $\mathds{Q}$-Cartier, and $g^+$-nef. 
\end{defin}
Pictorially, the result after the flop also gives rise to a morphism to $\Spec\scrR$, as follows.
\[
\begin{tikzpicture}
\node (X1) at (5,0) {$\scrX$};
\node (X2) at (8,0) {$\scrX^+$};
\node (Xcon) at (6.5,-1) {$\scrX_{\con}$};
\node (R) at (6.5,-2.5) {$\Spec \scrR$};
\draw[->] (X1) -- node[above, pos=0.6] {$\scriptstyle g$} (Xcon);
\draw[->] (X2) -- node[above,pos=0.6] {$\scriptstyle g^+$} (Xcon);
\draw[->] (Xcon) -- node[right] {$\scriptstyle $} (R);
\draw[->] (X1) -- node[below left] {$\scriptstyle f$} (R);
\draw[->] (X2) -- node[below right] {$\scriptstyle f'$} (R);
\end{tikzpicture}
\]
The idea in \ref{flopsdefin} is, roughly speaking, to use intersection theory to ensure that the same curve has not just been inserted straight back in, hence avoiding the identity operation. There are many other equivalent definitions of flops in the literature, see e.g.\ \cite{KollarFlops}.  

\begin{example}
Examples~\ref{toric ex}\eqref{toric ex 1} and \ref{nontoric ex}\eqref{nontoric ex 1} are examples of flops.  Example~\ref{nontoric ex}\eqref{nontoric ex 2} is not: it is a divisor-to-curve contraction, as are the three morphisms in Example~\ref{toric ex}\eqref{toric ex 2}.  However, there, the top curve in the left hand picture flops into the surface, giving the middle resolution.   Similarly, by symmetry the bottom curve in the right picture also flops into the surface.  These contractions, and flops, are is summarised as follows.
\[
\begin{array}{c}
\begin{tikzpicture}
\node at (-3,0) {\begin{tikzpicture}[scale=0.6]
\coordinate (T1) at (1.7,2.5);
\coordinate (T) at (1.9,1.8);
\coordinate (TM) at (1.92,1.7);
\coordinate (B) at (2.1,0.9);
\foreach \y in {0,0.1,0.2,...,1}{ 
\draw[line width=0.5pt,red!50] ($(T)+(\y,0)$) to [bend left=25] ($(B)+(\y,0)$);
\draw[line width=0.5pt,red!50] ($(T)+(-\y,0)$) to [bend left=25] ($(B)+(-\y,0)$);;}
\draw[red,line width=1pt] (T1) to [bend left=25] (TM);
\draw[line width=1pt,red] (T) to [bend left=25] (B);
\draw[color=blue!60!black,rounded corners=6pt,line width=0.75pt] (0.5,0) -- (1.5,0.3)-- (3.6,0) -- (4.3,1.5)-- (4,3.2) -- (2.5,2.7) -- (0.2,3) -- (-0.2,2)-- cycle;
\end{tikzpicture}};
\node at (0,0) {\begin{tikzpicture}[scale=0.6]
\coordinate (T) at (1.9,2);
\coordinate (TM) at (2.12-0.05,1.5-0.1);
\coordinate (BM) at (2.12-0.09,1.5+0.12);
\coordinate (B) at (2.1,1);
\foreach \y in {0.1,0.2,...,1}{ 
\draw[line width=0.5pt,red!50] ($(T)+(\y,0)+(0.02,0)$) to [bend left=25] ($(B)+(\y,0)+(0.02,0)$);
\draw[line width=0.5pt,red!50] ($(T)+(-\y,0)+(-0.02,0)$) to [bend left=25] ($(B)+(-\y,0)+(-0.02,0)$);;}
\draw[red,line width=1pt] (T) to [bend left=25] (TM);
\draw[red,line width=1pt] (BM) to [bend left=25] (B);
\draw[color=blue!60!black,rounded corners=6pt,line width=0.75pt] (0.5,0) -- (1.5,0.3)-- (3.6,0) -- (4.3,1.5)-- (4,3.2) -- (2.5,2.7) -- (0.2,3) -- (-0.2,2)-- cycle;
\end{tikzpicture}};
\node at (3,0) {\begin{tikzpicture}[scale=0.6]
\coordinate (T) at (1.9,2.2);
\coordinate (B) at (2.1,1.3);
\coordinate (BM) at (2.07,1.45);
\coordinate (B1) at (2.3,0.65);
\foreach \y in {0,0.1,0.2,...,1}{ 
\draw[line width=0.5pt,red!50] ($(T)+(\y,0)$) to [bend left=25] ($(B)+(\y,0)$);
\draw[line width=0.5pt,red!50] ($(T)+(-\y,0)$) to [bend left=25] ($(B)+(-\y,0)$);;}
\draw[red,line width=1pt] (BM) to [bend left=25] (B1);
\draw[line width=1pt,red] (T) to [bend left=25] (B);
\draw[color=blue!60!black,rounded corners=6pt,line width=0.75pt] (0.5,0) -- (1.5,0.3)-- (3.6,0) -- (4.3,1.5)-- (4,3.2) -- (2.5,2.7) -- (0.2,3) -- (-0.2,2)-- cycle;
\end{tikzpicture}};
\node at (-1.5,-2.25) {\begin{tikzpicture}[scale=0.6]
\coordinate (T) at (1.9,2);
\coordinate (B) at (2.1,1);
\foreach \y in {0,0.1,0.2,...,1}{ 
\draw[line width=0.5pt,red!50] ($(T)+(\y,0)+(0.02,0)$) to [bend left=25] ($(B)+(\y,0)+(0.02,0)$);
\draw[line width=0.5pt,red!50] ($(T)+(-\y,0)+(-0.02,0)$) to [bend left=25] ($(B)+(-\y,0)+(-0.02,0)$);;}
\draw[red,line width=1pt] (T) to [bend left=25] (B);
\filldraw[blue!60!black] (T) circle (2pt); 
\draw[color=blue!60!black,rounded corners=6pt,line width=0.75pt] (0.5,0) -- (1.5,0.3)-- (3.6,0) -- (4.3,1.5)-- (4,3.2) -- (2.5,2.7) -- (0.2,3) -- (-0.2,2)-- cycle;
\end{tikzpicture}};
\node at (1.5,-2.25) {\begin{tikzpicture}[scale=0.6]
\coordinate (T) at (1.9,2);
\coordinate (B) at (2.1,1);
\foreach \y in {0,0.1,0.2,...,1}{ 
\draw[line width=0.5pt,red!50] ($(T)+(\y,0)+(0.02,0)$) to [bend left=25] ($(B)+(\y,0)+(0.02,0)$);
\draw[line width=0.5pt,red!50] ($(T)+(-\y,0)+(-0.02,0)$) to [bend left=25] ($(B)+(-\y,0)+(-0.02,0)$);;}
\draw[red,line width=1pt] (T) to [bend left=25] (B);
\filldraw[blue!60!black] (B) circle (2pt); 
\draw[color=blue!60!black,rounded corners=6pt,line width=0.75pt] (0.5,0) -- (1.5,0.3)-- (3.6,0) -- (4.3,1.5)-- (4,3.2) -- (2.5,2.7) -- (0.2,3) -- (-0.2,2)-- cycle;
\end{tikzpicture}};
\draw[->, color=blue!60!black] (-2.25,-1) -- (-1.75,-1.4);
\draw[->, color=blue!60!black] (-0.75,-1) -- (-1.25,-1.4);
\draw[->, color=blue!60!black] (0.75,-1) -- (1.25,-1.4);
\draw[->, color=blue!60!black] (2.25,-1) -- (1.75,-1.4);
\end{tikzpicture}
\end{array}
\]
\end{example}

\subsection{Terminal singularities and all that}\label{terminal}
In the wider perspective of the MMP, cDV singularities play various roles.  From the viewpoint of this survey, the two key concepts are that of Gorenstein canonical and Gorenstein terminal singularities. 

Both concepts have definitions, found in \cite{YPG}, and we do not seek to repeat them here. It turns out, via Grothendieck duality, that Gorenstein canonical is equivalent to being Gorenstein and rational.  That is to say, $R$ is Gorenstein canonical if and only if $R$ is Gorenstein, and furthermore there exists a resolution $X\to\Spec R$ for which $\RDerived f_*\scrO_X=\scrO_R$.  It is equivalent to insist that all resolutions have this property.  The property $\RDerived f_*\scrO_X=\scrO_R$ really is key, from a noncommutative perspective, since it is needed for tilting to exist. It is from tilting that the noncommutative aspects of the theory are born; see \S\ref{NCCR section}.

Perhaps the next most important characterisation is Reid's original observation that \emph{isolated} cDV singularities are precisely the Gorenstein terminal singularities, or `terminal of index one'. This gives rise to the following sequence of equalities and inclusions
\[
\{ \textnormal{Gorenstein terminal}\}=
\{ \textnormal{isolated cDV}\}\subset
\{ \textnormal{all cDV}\} \subset
\{ \textnormal{Gorenstein canonical}\} .
\]
Two points stand out.  First, it is the Gorenstein terminal characterisation of isolated cDV singularities that links cDV singularities firmly to the theory of flops.  Second, it is almost impossible to convey in words how large the class of Gorenstein canonical singularities is \cite{ReidCanonical}.   The cDV singularities are an extremely small subset, and are by no means representative of the general complexity.

\medskip
The philosophy, from the homological perspective, is that there should be no difference between the setting of crepant partial resolutions $\scrX\to\Spec\scrR$ of a cDV $\scrR$ and the specific case when $\scrX$ is smooth.  In the flops case, this translates into the observation that unified statements should be expected for all flopping contractions $\scrX\to\Spec\scrR$, where $\scrX$ has at worst Gorenstein terminal singularities.

\subsection{The General Elephant}
The passage between $3$-fold cDV singularities and that of their surface counterparts is provided by Reid's general elephant.  The key result \cite{Pagoda} is that for \emph{generic} $g\in\scrR$, for any crepant partial resolution $\scrX\to\Spec\scrR$ the pullback diagram
\begin{equation}
\begin{array}{c}
\begin{tikzpicture}
\node (A) at (0,0) {$\scrY$};
\node (a) at (2,0) {$\scrX$};
\node (B) at (0,-1.25) {$\Spec\scrR/g$};
\node (b) at (2,-1.25) {$\Spec\scrR$};
\draw[->] (A)--(a);
\draw[->] (A)--(B);
\draw[->] (B)--(b);
\draw[->] (a)--(b);
\end{tikzpicture}
\end{array}\label{general elephant}
\end{equation}
results in an $\scrR/g$ which is an ADE surface singularity, \emph{and} $\scrY\to\Spec\scrR/g$ is a crepant partial resolution.  An often misunderstood point, which complicates the combinatorics below, is that $\scrY\to\Spec\scrR/g$ need not be the minimal resolution.

\begin{example}\label{cut ex}
We illustrate the above with two examples.
\begin{enumerate}
\item\label{cut ex 1} Consider $\scrR=\mathbb{C}\llsq u,v,x,y\rrsq/(uv-(x^2+y^3)x)$, and its minimal models described in \ref{nontoric ex}\eqref{nontoric ex 1}.  Under generic pullback, the flop picture in \ref{nontoric ex}\eqref{nontoric ex 1} results in the two different singular partial resolutions of the $A_2$ surface singularity, which are dominated by the minimal resolution as follows.
\[
\begin{tikzpicture}
\node at (-1.25,0) {$\begin{tikzpicture}[scale=0.5]
\coordinate (T) at (1.9,2);
\coordinate (B) at (2.1,1);
\draw[red,line width=1pt] (T) to [bend left=25] (B);
\filldraw[blue!60!black] (B) circle (2pt); 
\draw[color=blue!60!black,rounded corners=5pt,line width=1pt] (0.5,0) -- (1.5,0.3)-- (3.6,0) -- (4.3,1.5)-- (4,3.2) -- (2.5,2.7) -- (0.2,3) -- (-0.2,2)-- cycle;
\draw[densely dotted] (1.25,2.5) -- (3.5,2.5) -- (3,0.5) -- (0.75,0.5) -- cycle; 
\end{tikzpicture}$};
\node at (1.25,0) {$\begin{tikzpicture}[scale=0.5]
\coordinate (T) at (2.1,2);
\coordinate (B) at (1.9,1);
\draw[red,line width=1pt] (T) to [bend right=25] (B);
\filldraw[blue!60!black] (T) circle (2pt); 
\draw[color=blue!60!black,rounded corners=5pt,line width=1pt] (0.5,0) -- (1.5,0.3)-- (3.6,0) -- (4.3,1.5)-- (4,3.2) -- (2.5,2.7) -- (0.2,3) -- (-0.2,2)-- cycle;
\draw[densely dotted] (1.25,2.5) -- (3.5,2.5) -- (3,0.5) -- (0.75,0.5) -- cycle; 
\end{tikzpicture}$};
\node at (0,-2) {\begin{tikzpicture}[scale=0.5]
\filldraw [red] (2.1,0.75) circle (1pt);
\draw[color=blue!60!black,rounded corners=5pt,line width=1pt] (0.5,0) -- (1.5,0.15)-- (3.6,0) -- (4.3,0.75)-- (4,1.6) -- (2.5,1.35) -- (0.2,1.5) -- (-0.2,1)-- cycle;
\draw[densely dotted] (1.25,1) -- (3.75,1) -- (3.25,0.5) -- (0.75,0.5) -- cycle; 
\end{tikzpicture}};
\draw[->, color=blue!60!black] (-0.75,-1) -- (-0.25,-1.5);
\draw[->, color=blue!60!black] (0.75,-1) -- (0.25,-1.5);
\node at (-7.25,0) {$\begin{tikzpicture}[scale=0.5]
\coordinate (T) at (1.9,2);
\coordinate (B) at (2.1,1);
\draw[red,line width=1pt] (T) to [bend left=25] (B);
\filldraw[blue!60!black] (B) circle (2pt); 
\draw[draw=none,color=blue!60!black,rounded corners=5pt,line width=1pt] (0.5,0) -- (1.5,0.3)-- (3.6,0) -- (4.3,1.5)-- (4,3.2) -- (2.5,2.7) -- (0.2,3) -- (-0.2,2)-- cycle;
\draw[color=blue!60!black,,line width=1pt] (1.25,2.5) -- (3.5,2.5) -- (3,0.5) -- (0.75,0.5) -- cycle; 
\end{tikzpicture}$};
\node at (-4.75,0) {$\begin{tikzpicture}[scale=0.5]
\coordinate (T) at (2.1,2);
\coordinate (B) at (1.9,1);
\draw[red,line width=1pt] (T) to [bend right=25] (B);
\filldraw[blue!60!black] (T) circle (2pt); 
\draw[draw=none,color=blue!60!black,rounded corners=5pt,line width=1pt] (0.5,0) -- (1.5,0.3)-- (3.6,0) -- (4.3,1.5)-- (4,3.2) -- (2.5,2.7) -- (0.2,3) -- (-0.2,2)-- cycle;
\draw[color=blue!60!black,,line width=1pt]  (1.25,2.5) -- (3.5,2.5) -- (3,0.5) -- (0.75,0.5) -- cycle; 
\end{tikzpicture}$};
\node at (-6,-2) {\begin{tikzpicture}[scale=0.5]
\filldraw [red] (2.1,0.75) circle (1pt);
\draw[draw=none,color=blue!60!black,rounded corners=5pt,line width=1pt] (0.5,0) -- (1.5,0.15)-- (3.6,0) -- (4.3,0.75)-- (4,1.6) -- (2.5,1.35) -- (0.2,1.5) -- (-0.2,1)-- cycle;
\draw[color=blue!60!black,line width=1pt]  (1.25,1) -- (3.75,1) -- (3.25,0.5) -- (0.75,0.5) -- cycle; 
\end{tikzpicture}};
\node at (-6,2) {$\begin{tikzpicture}[scale=0.5]
\coordinate (T1) at (1.9,2.5);
\coordinate (T) at (2.25,1.85);
\coordinate (TM) at (2.2,1.7);
\coordinate (B) at (2,0.9);
\draw[red,line width=1pt] (T1) to [bend left=25] (TM);
\draw[line width=1pt,red] (T) to [bend right=25] (B);
\draw[draw=none,color=blue!60!black,rounded corners=5pt,line width=1pt] (0.5,0) -- (1.5,0.3)-- (3.6,0) -- (4.3,1.5)-- (4,3.2) -- (2.5,2.7) -- (0.2,3) -- (-0.2,2)-- cycle;
\draw[color=blue!60!black,line width=1pt]  (1.25,2.75) -- (3.5,2.75) -- (3,0.5) -- (0.75,0.5) -- cycle; 
\end{tikzpicture}$};
\draw[->, color=blue!60!black] (-6.75,-1) -- (-6.25,-1.5);
\draw[->, color=blue!60!black] (-5.25,-1) -- (-5.75,-1.5);
\draw[->, color=blue!60!black] (-6.25,1.25) -- (-6.75,0.75);
\draw[->, color=blue!60!black] (-5.75,1.25) -- (-5.25,0.75);
\draw[color=blue!60!black,right hook->,line width=0.75pt] (-3.5,-1.25) -- (-2.75,-1.25);
\node at (-8,-2) {$\Spec(\scrR/g)$};
\node at (2,-2) {$\Spec\scrR$};
\end{tikzpicture}
\]
This example illustrates two points: (1) $\scrY$ need not be smooth, and (2) the surface partial resolution obtained by slicing varies depending on the minimal model.  How to obtain one partial resolution from the other is controlled by the combinatorial \emph{wall crossing rule}, described in \S\ref{wall crossing section}.
\item\label{cut ex 2} (Laufer flop) Consider $\scrR=\mathbb{C}\llsq u,v,x,y\rrsq/(u^2+v^2y-(x^2+y^3)x)$, which is isolated $cD_4$.  Blowing up the reflexive ideal $(x^2+y^3, vx+uy,ux-vy^2)$ gives a crepant resolution $\scrX\to\Spec\scrR$, where the unique singular point has been replaced by a single curve.  Under generic slice, $\scrR/g$ is $D_4$, and since there is only one curve, the pullback $\scrY\to\Spec\scrR/g$ cannot be the minimal resolution.
\end{enumerate}
\end{example}

\subsection{Ample and Movable Cones}
Two key pieces of birational geometry data associated to $f\colon \scrX\to\Spec\scrR$ are the ample and movable cones.  We partially recap the definitions here, mainly to set notation, then return to them in more detail in \S\ref{Pinkham sect}.

Consider $\mathrm{C}\colonequals f^{-1}(p)$, the fibre above the unique closed point. Equipping $\mathrm{C}$ with its reduced scheme structure, we obtain $\mathrm{C}^{\rm red}=\bigcup_{i=1}^n\mathrm{C}_i$ with each $\mathrm{C}_i\cong\mathbb{P}^1$. Since $\scrR$ is complete there exist $D_i$ such that $D_i\cdot \mathrm{C}_j=\updelta_{ij}$. Two divisors $D, E$ on $\scrX$ are called \emph{numerically equivalent} over $\scrR$ if $D\cdot\mathrm{C}_i = E \cdot \mathrm{C}_i$ for all $i$.  With these notions set,
\[
\mathrm{N}^1(\scrX/\scrR)=\frac{\{\textnormal{divisors on }\scrX\}}{\textnormal{numerical equivalence over }\scrR}\otimes_{\mathbb{Z}}\mathbb{R}
\cong
\bigoplus_{i=1}^n\mathbb{R}[D_i]
\]
is called the \emph{group of divisor classes}.  The \emph{ample cone} is defined to be
\[
\mathrm{A}^1(\scrX/\scrR)=\bigoplus_{i=1}^n\mathbb{R}_{> 0}[D_i].
\]
If $g\colon \scrW\to\Spec\scrR$ can obtained from the fixed $\scrX\to \Spec\scrR$ by iteratively flopping curves, consider the induced birational map $\upphi=f^{-1}\circ g\colon \scrW\birational\scrX$. Being an isomorphism in codimension one, there is an induced isomorphism $\upphi_*\colon  \mathrm{N}^1(\scrW/\scrR)\to\mathrm{N}^1(\scrX/\scrR)$, called the strict transform.

In the flops case, equivalently when $\scrR$ is isolated, by \cite[2.3]{KawamataCone} we have
\begin{equation}
\mathrm{N}^1(\scrX/\scrR)=\bigcup_{\scrW}\,\upphi_*\overline{\mathrm{A}(\scrW/\scrR)}\label{KawaCone}
\end{equation}
where the union is over all $\scrW\to\Spec\scrR$ obtained from $\scrX\to\Spec\scrR$ by iteratively flopping curves.  In the special case when $\scrX$ is a minimal model, the union is equivalently over the set of all minimal models.

\begin{remark} 
(See \cite[3.4]{TodaResPub}) In the flops case, $f \colon\scrX\to\Spec\scrR$ is an isomorphism in codimension one, so the $f$-effective $f$-movable cone and $f$-effective $f$-nef cone defined in \cite[1.1]{KawamataCone} coincide with $\mathrm{N}^1(\scrX/\scrR )$ and $\mathrm{A}(\scrX/\scrR)$ respectively.
\end{remark}


\section{The Combinatorics of cDV Singularities}\label{combo section}

Much of the birational geometry, and homological algebra, of cDV singularities can be encoded in certain combinatorial objects called \emph{intersection} (or restriction) arrangements.  These are constructed from ADE data, but need not themselves be ADE;  we motivate their construction geometrically in \S\ref{Pinkham sect}.  From these arrangements, and the wall crossing rule in \S\ref{wall crossing section}, many geometric properties can be extracted.

\subsection{Intersection Arrangement Generalities}\label{int arr section}
Consider an ADE Dynkin diagram $\Updelta$, with associated Coxeter group, denoted $\WDelt$.  Let $V$ be the $\mathbb{R}$-vector space with basis $\{\upalpha_i\mid i\in \Updelta\}$, and set $\Uptheta\colonequals V^*$ to be the dual space of $V$, which has basis $\{\upalpha_i^* \mid i\in \Updelta\}$.   The Weyl group acts on $V^*$, and the \emph{Tits cone} $\Cone{\Updelta}$ is defined to be
\[
\Cone{\Updelta}\colonequals \bigcup_{x\in \WDelt}x(\overline{C}),
\]
where $\overline{C}\colonequals \{\upvartheta\in\Uptheta\mid \upvartheta_i\geq 0 \mbox{ for all }i\in\Updelta\}$.

For any choice $\scrJ$ of a subset of vertices of $\Updelta$, consider the vector space
\[
\Uptheta_{\scrJ}\colonequals \{ \upvartheta\in\Uptheta\mid \upvartheta_i=0\mbox{ if }i\in \scrJ\},
\]
which has as basis $\{\upalpha_i^*\mid i\notin\scrJ\}$.   The main object of our study, called the $\scrJ$-cone, is the intersection
\[
\Cone{\Updelta,\scrJ}\colonequals \Cone{\Updelta}\cap \Uptheta_{\scrJ}.
\]
The hyperplanes of $\Cone{\Updelta}$ intersect the subspace $\Cone{\Updelta,\scrJ}$ to form a new hyperplane arrangement, called the \emph{intersection arrangement} $\scrH_\scrJ$.  This will often be simplified to $\scrH$, whenever the $\scrJ$ is clear.

In order to describe $\scrH_\scrJ$, it is convenient to label its chambers.  For $\scrJ\subseteq\Updelta$, let $\cham(\Updelta,\scrJ)$ be the set of those pairs $(x,\scrI)$, where $x$ is an element of the Weyl group $\WDelt$,  $\scrI$ is a subset of $\Updelta$, satisfying the two properties that $\ell(x)=\min\{\ell(y)\mid y\in x\Wkern{\scrI}\}$, and $\Wkern{\scrJ} x=x\Wkern{\scrI}$.

When $\Updelta$ is a Coxeter graph, it turns out that there is a bijection from the set $\cham(\Updelta,\scrJ)$ to the set of chambers in $\scrH_\scrJ$, given by
\begin{equation}
\cham(\Updelta,\scrJ)\ni(x,\scrI)\mapsto x(C_{\scrI})\label{Cox labels}
\end{equation}
This allows us to label the chambers by Coxeter data, even when the arrangement $\scrH_{\scrJ}$ itself need not be Coxeter.  This labelling becomes important geometrically in \S\ref{wall crossing section}, since it describes how to obtain all partial resolutions in the generic slice from any fixed one.  We remark that in the classical case, namely when $\scrJ=\emptyset$, the above \eqref{Cox labels} reduces to the well-known fact that elements of the Weyl group label chambers.

\medskip
Consider next the extended Dynkin case $\Updelta_{\aff}$.  Given a subset $\scrJ$ of vertices of the Dynkin diagram $\Updelta$, we can also consider $\scrJ$ as a subset of $\Updelta_{\aff}$.  Applying the above to $\scrJ\subset\Updelta_{\aff}$, we then call
\[
\Cone{\scrJ_{\aff}}\colonequals \Cone{\Updelta_{\aff},\scrJ}\subseteq\mathbb{R}^{|\Updelta_{\aff}|-|\scrJ|}
\]
the \emph{$\scrJ$-affine Tits cone}.   As for usual Tits cones in affine settings, there is redundancy as $\Cone{\scrJ_{\aff}}$ does not fill $\mathbb{R}^{|\Updelta_{\aff}|-|\scrJ|}$.  Consequently, the \emph{level} is defined to be
\begin{equation}
\Level(\scrJ)\colonequals \{ \upvartheta\in \Cone{\Updelta_{\aff},\scrJ}\mid \sum_{k\notin\scrJ} \updelta_k\upvartheta_k=1\}.\label{level}
\end{equation}
where $(\updelta_i)_{i\in\Delt_{\aff}}$ is the imaginary root. Thus for $\scrJ\subseteq\Delt$, the level $\Level(\scrJ_{\aff})$ inside $\Cone{\Updelta_{\aff},\scrJ}$  is an infinite hyperplane arrangement in $\mathbb{R}^{|\Delt|-|\scrJ|}$.  We write $\scrH^{\aff}$ for this infinite hyperplane arrangement. The chambers of $\Cone{\scrJ_{\aff}}$ partition $\Level(\scrJ_{\aff})$ into \emph{alcoves}.   The easiest example is given in \cite[2.9]{IW9}.

\medskip
The upshot of the above is that for any choice $\scrJ$ inside an ADE Dynkin diagram $\Delt$, there is an associated finite hyperplane arrangement $\scrH$, and an infinite arrangement $\scrH^{\aff}$, both living inside $\mathbb{R}^{|\Delt|-|\scrJ|}$

\begin{example}
If $\Updelta=E_6$ and $\scrJ=\Esix{R}{R}{W}{W}{R}{R}$,  where the nodes shaded red are the nodes in $\scrJ$, then $\scrH$ and $\scrH^{\aff}$ are the arrangements in $\mathbb{R}^{|\Delt|-|\scrJ|}=\mathbb{R}^2$ illustrated below.
\[
\begin{array}{ccc}
\begin{array}{c}
\begin{tikzpicture}[scale=0.6]
\clip (-3.3,-3.3) rectangle (2.3,2.3);
{\foreach \y in {-1,0,1}\foreach \x in {-1,0,1}
{
\node at (3*\y,3*\x) 
{\begin{tikzpicture}[scale=0.6]
\draw[draw=none] (0,0) -- (0,3) -- (3,3) -- (3,0) --cycle;
\draw[draw=none] (0,1) -- (3,1);
\draw[draw=none] (0,2) -- (3,2);
\draw[draw=none] (1,0) -- (1,3);
\draw[draw=none] (2,0) -- (2,3);
\coordinate (l) at (0,1.5);
\coordinate (u) at (1.5,3);
\coordinate (d) at (1.5,0);
\coordinate (r) at (3,1.5);
\draw[draw=none] (0,3) -- (d);
\draw[draw=none] (0,3) -- (r);
\draw[draw=none] (3,0) -- (u);
\draw[draw=none] (3,0) -- (l);
\draw[draw=none] (l) -- (u);
\draw[draw=none] (0,0) -- (3,3);
\draw[draw=none] (d) -- (r);
\end{tikzpicture}};}}
\draw[green!60!black,thick] (-0.5,5)--(-0.5,-5);
\draw[green!60!black,thick] (-4.5,-0.5)--(4,-0.5);
\draw[green!60!black,thick] (3.5,3.5)--(-4.5,-4.5);
\draw[green!60!black,thick] (-3.25,5)--(1.75,-5);
\draw[green!60!black,thick] (-4.5,1.5)--(4,-2.75);
\end{tikzpicture}
\end{array}
&&
\begin{array}{c}
\begin{tikzpicture}[scale=0.6]
\clip (-3.3,-3.3) rectangle (2.3,2.3);
{\foreach \y in {-1,0,1}\foreach \x in {-1,0,1}
{
\node at (3*\y,3*\x) 
{\begin{tikzpicture}[scale=0.6]
\draw (0,0) -- (0,3) -- (3,3) -- (3,0) --cycle;
\draw (0,1) -- (3,1);
\draw (0,2) -- (3,2);
\draw (1,0) -- (1,3);
\draw (2,0) -- (2,3);
\coordinate (l) at (0,1.5);
\coordinate (u) at (1.5,3);
\coordinate (d) at (1.5,0);
\coordinate (r) at (3,1.5);
\draw (0,3) -- (d);
\draw (0,3) -- (r);
\draw (3,0) -- (u);
\draw (3,0) -- (l);
\draw (l) -- (u);
\draw (0,0) -- (3,3);
\draw (d) -- (r);
\end{tikzpicture}};}}
\draw[green!60!black,thick] (-0.5,5)--(-0.5,-5);
\draw[green!60!black,thick] (-4.5,-0.5)--(4,-0.5);
\draw[green!60!black,thick] (3.5,3.5)--(-4.5,-4.5);
\draw[green!60!black,thick] (-3.25,5)--(1.75,-5);
\draw[green!60!black,thick] (-4.5,1.5)--(4,-2.75);
\end{tikzpicture}
\end{array}
\end{array}
\]
\end{example}

\begin{example}
If $\Updelta=E_8$ and $\scrJ=\Eeight{R}{R}{W}{R}{R}{W}{R}{W}$, where again the nodes shaded red are the nodes in $\scrJ$, then the corresponding finite arrangement $\scrH$ in $\mathbb{R}^{|\Delt|-|\scrJ|}=\mathbb{R}^3$ is illustrated in \eqref{192 picture}, and has 192 chambers.  
\end{example}

\subsection{The Wall Crossing Rule}\label{wall crossing section}
Combinatorially, if in \eqref{Cox labels} $(x,\scrI)$ and $(y,\scrI')$ label adjacent chambers, then it is possible to describe one from the other combinatorially, via a \emph{wall crossing rule}.  Here we skip over the rule for obtaining $y$ from $x$ (and vice versa), and instead focus on how to obtain $\scrI'$ from $\scrI$.  Full details can be found in \cite{IW9}.

To describe this, recall first that ADE Dynkin diagrams carry a canonical involution $\upiota$, usually referred to as the \emph{Dynkin involution}.  In Type $A_n$ it acts as a reflection in the centre point of the chain, which may or may not be a vertex, namely
\[
\begin{tikzpicture}[>=stealth,scale=0.45]
 \node (1) at (1,0) [DB] {};
 \node (2) at (2,0) [DB] {};
 \node (3) at (4,0) [DB] {};
 \node (4) at (5,0)[DB] {};
\draw [-] (1) -- (2);
\draw [-] (2) -- (2.6,0);
\draw [-] (3.4,0) -- (3);
\draw [-] (3) -- (4);
\draw[densely dotted] (3,0.5)--(3,-0.5);
\draw[<->, bend left,densely dotted] (2.5,0.5) to (3.5,0.5);
\end{tikzpicture}
\]
When $n\geq 4$ and $n$ is odd, the involution $\upiota$ acts on $D_{n}$ by permuting the left hand branches, and fixing all other vertices:
\[
\begin{tikzpicture}[>=stealth,scale=0.45]
 \node (0) at (0,0) [DB] {};
 \node (1) at (1,0) [DB] {};
 \node (1b) at (1,0.75) [DB] {};
 \node (2) at (2,0) [DB] {};
 \node (3) at (4,0) [DB] {};
 \node (4) at (5,0)[DB] {};
 \node at (3,0) {$\cdots$};
\draw [-] (0) -- (1);
\draw [-] (1) -- (2);
\draw [-] (2) -- (2.6,0);
\draw [-] (3.4,0) -- (3);
\draw [-] (3) -- (4);
\draw [-] (1) -- (1b);
\draw[<->, bend left,densely dotted] (0) to (1b);
\end{tikzpicture}
\]
In type $E_6$ the involution $\upiota$ acts as a reflection, and in all other cases, $\upiota$ is the identity.

\medskip
Now, consider an arbitrary chamber, say with label $(x,\scrI)$.  This chamber, like all chambers, has walls which correspond to the nodes not in $\scrI$.  To cross any one of these walls is equivalent to choosing a node not in $\scrI$, and applying the combinatorial rule below. In our conventions, these are the nodes \emph{not} shaded red.

Hence, to cross one of these walls, choose a node not in $\scrI$, i.e.\ a non-red node.  Temporarily delete \emph{all other} non-red nodes, apply the Dynkin involution $\upiota$, then put back in the deleted vertices.  This is best illustrated via an example.

\begin{example}
Suppose $\scrI$ is the following, where we have circled the non-red node at which we wish to wall-cross.
\[
\begin{tikzpicture}[scale=0.45]
 \node (m1) at (-1,0) [DB] {};
 \node (0) at (0,0) [DB] {};
 \node (1) at (1,0) [DW] {};
 \node (1b) at (1,0.75) [DB] {};
 \node (2) at (2,0) [DW] {};
 \node (3) at (3,0) [DB] {};
 \node (4) at (4,0) [DB] {};
  \node (5) at (5,0) [DW] {};
\draw[densely dotted, circle] (4) circle (14pt);
\draw [-] (m1) -- (0);
\draw [-] (0) -- (1);
\draw [-] (1) -- (2);
\draw [-] (2) -- (3);
\draw [-] (3) -- (4);
\draw [-] (4) -- (5);
\draw [-] (1) -- (1b);
\end{tikzpicture}
\]
Deleting all other non-red nodes we are left with the following, on which $\upiota$ acts as indicated. 
\[
\begin{tikzpicture}[scale=0.45,>=stealth]
 \node (m1) at (-1,0) [] {};
 \node (0) at (0,0) [] {};
 \node (1) at (1,0) [DW] {};
 \node (2) at (2,0) [DW] {};
 \node (3) at (3,0) [] {};
 \node (4) at (4,0) [DB] {};
  \node (5) at (5,0) [DW] {};
 \draw [-] (1) -- (2);
 \draw [-] (4) -- (5);
\draw[<->] (1,-0.7) to (2,-0.7);
\draw[<->] (4,-0.7) to (5,-0.7);
  \end{tikzpicture}
\]
Applying the Dynkin involution then moves the non-red node to the right:
\[
\begin{tikzpicture}[scale=0.45,>=stealth]
 \node (m1) at (-1,0) [] {};
 \node (0) at (0,0) [] {};
 \node (1) at (1,0) [DW] {};
 \node (2) at (2,0) [DW] {};
 \node (3) at (3,0) [] {};
 \node (4) at (4,0) [DW] {};
  \node (5) at (5,0) [DB] {};
 \draw [-] (1) -- (2);
 \draw [-] (4) -- (5);
  \end{tikzpicture}
\]
The final step is to add back in the temporarily deleted vertices, to obtain the following.
\[
\begin{tikzpicture}[scale=0.45]
 \node (m1) at (-1,0) [DB] {};
 \node (0) at (0,0) [DB] {};
 \node (1) at (1,0) [DW] {};
 \node (1b) at (1,0.75) [DB] {};
 \node (2) at (2,0) [DW] {};
 \node (3) at (3,0) [DB] {};
 \node (4) at (4,0) [DW] {};
  \node (5) at (5,0) [DB] {};
\draw [-] (m1) -- (0);
\draw [-] (0) -- (1);
\draw [-] (1) -- (2);
\draw [-] (2) -- (3);
\draw [-] (3) -- (4);
\draw [-] (4) -- (5);
\draw [-] (1) -- (1b);
\end{tikzpicture}
\]
\end{example}

It is a somewhat remarkable fact that the above combinatorial rule is how generic slices of different minimal models are related \cite{HomMMP, IW9}.  Indeed, if $\scrX\to\Spec\scrR$ is some crepant partial resolution, such that the slice is controlled by $\scrI\subset\Delt$, then the flops of $\scrX$ at  irreducible curves slice precisely to the different wall crossings of $\scrI$.

\subsection{Pinkham}\label{Pinkham sect}
Given a crepant partial resolution $\scrX\to\Spec\scrR$, where $\scrR$ is cDV, the relevant combinatorics are controlled by the general elephant \eqref{general elephant}.  Indeed, taking the pullback with respect to a generic $g$ gives rise to a crepant partial resolution $\scrY\to\Spec\scrR/g$, where $\scrR/g$ is ADE.  Thus $\scrY$ is dominated by the minimal resolution
\[
\scrY_{\mathrm{min}}\xrightarrow{h}\scrY\to\Spec\scrR/g.
\]
Consequently, to describe $\scrY$, it suffices to describe which curves are contracted under the morphism $h$.  This is combinatorial.  Since the dual graph of $\scrY_{\mathrm{min}}\to\Spec\scrR/g$ is an ADE Dynkin diagram $\Delt$, to say which curves are contracted by $h$ is equivalent to prescribing this subset of $\Delt$.  We write $\scrJ$ for the subset of $\Delt$ which correspond to the curves contracted by $h$, and thus via the previous subsection obtain  hyperplane arrangements $\scrH$ and $\scrH^{\aff}$ in $\mathbb{R}^{|\Delt|-|\scrJ|}$. We remark that $|\Delt|-|\scrJ|$ is the number of curves that \emph{survive} under $h$, which is equal to the number curves in $\scrX\to\Spec\scrR$ above the origin. 

\medskip
In the case $\scrX\to\Spec\scrR$ is a flopping contraction, so necessarily $\scrR$ is isolated, it is a theorem of Pinkham \cite[Theorem 3]{Pinkham} that the movable cone is the hyperplane arrangement $\scrH$.  This is a more precise description than that given by \eqref{KawaCone}, since it describes the cones in terms of the intersection arrangement combinatorics.  

This Pinkham result was recovered and generalised in \cite[\S6]{HomMMP}, where it was proved that for an arbitrary cDV $\scrR$,  in the skeleton of $\scrH$ there exists a connected path, containing the chamber labelled $(1,\scrJ)$, where every minimal model of $\scrR$ appears precisely once, and furthermore each wall crossing in this path corresponds to a flop.  Thus, in general there exists a connected chain inside $\scrH$ which describes the birational geometry.

\begin{example}\label{7.7HomMMPa}
This is \cite[7.7]{HomMMP}.  The $cD_4$ singularity $\scrR=\mathbb{C}\llsq u,v,y,z\rrsq /(u^2-v(x^2-4y^3))$ has precisely two crepant resolutions. These two minimal models of $\Spec\scrR$ can be viewed inside the finite hyperplane arrangement $\scrH$ as follows.
\[
\begin{array}{ccc}
\begin{array}{c}
\begin{tikzpicture}[scale=0.75]
\coordinate (A1) at (135:2cm);
\coordinate (A2) at (-45:2cm);
\coordinate (B1) at (153.435:2cm);
\coordinate (B2) at (-26.565:2cm);
\draw[red] (A1) -- (A2);
\draw[green!70!black] (B1) -- (B2);
\draw[-] (-2,0)--(2,0);
\draw[-] (0,-2)--(0,2);
\draw[densely dotted,gray] (0,0) circle (2cm);
\filldraw[blue] (45:1cm) circle (2pt);
\filldraw[blue] (112.5:1cm) circle (2pt);
\draw[blue] (45:1cm)--(112.5:1cm);
\end{tikzpicture}
\end{array}
\end{array}
\]
In this example, the corresponding infinite arrangement $\scrH^{\aff}$ is pictured in \cite[4.15]{IW9}.
\end{example}

\section{NC Resolutions and Variants}\label{NCCR section}

Noncommutative resolutions are the language in which we control iterations, through mutations and flops, for cDV singularities.  The major benefit is that whilst \emph{some} iterations come from birational geometry, most do not, and so noncommutative resolutions give a more holistic view on the derived category.  This extra information has geometric corollaries, mainly related to stability conditions and curve counting.

\subsection{Generalities}
Given $\scrR$ \textnormal{cDV} as before, $M\in\mod \scrR$ is called \emph{maximal Cohen-Macaulay} (=CM) provided
\[
\depth_\scrR M\colonequals\inf \{ i\geq 0\mid \Ext^i_\scrR(\scrR/\m,M)\neq 0 \}=\dim \scrR,
\]
and we write $\CM \scrR$ for the category of CM $\scrR$-modules.  Further, for $(-)^*\colonequals\Hom_\scrR(-,\scrR)$, $M\in\mod \scrR$ is called reflexive if the natural morphism $M\to M^{**}$ is an isomorphism, and we write $\refl \scrR$ for the category of reflexive $\scrR$-modules.

\begin{defin}\label{MMdefin}
We say $N\in\refl \scrR$ is a \emph{modifying module} if $\End_\scrR(N)\in\CM \scrR$, and we say that $N\in\refl \scrR$ is a \emph{maximal modifying (MM) module} if it is modifying and it is maximal with respect to this property;  equivalently,
\[
\add N=\{ X\in\refl \scrR\mid \End_{\scrR}(N \oplus X)\in\CM \scrR  \}.
\]
If $N$ is an \textnormal{MM} module, we call $\End_\scrR(N)$ a \emph{maximal modification algebra (=MMA).}  
\end{defin}

Of particular interest are members of $(\mathrm{MM}\,\scrR)\cap(\CM\scrR)$, namely those maximal modifying modules which are also CM. These are the maximal modifying \emph{generators} referred to below.  It will turn out that there are only finitely many.

The notion of a smooth noncommutative minimal model, called a noncommutative crepant resolution, is due to Van den Bergh \cite{VdBNCCR}.
\begin{defin}
A \emph{noncommutative crepant resolution (NCCR)} of $\scrR$ is $\Lambda\colonequals\End_\scrR(N)$ where $N\in\refl \scrR$ is such that $\Lambda\in\CM \scrR$ and $\gl\Lambda=3$. 
\end{defin}

It turns out that if there exists an NCCR $\End_\scrR(N)$, then $N$ is automatically MM, and further all MM modules give NCCRs.  In other words, if one noncommutative minimal model is smooth, they all are \cite[5.11]{IW3}.  Further, in the case that an NCCR exists, there is a more intrinsic characterisation of  $(\mathrm{MM}\,\scrR)\cap(\CM\scrR)$, namely
\[
N\in(\mathrm{MM}\,\scrR)\cap(\CM\scrR)
\iff
\add N=\{ X\in\CM \scrR \mid \Hom_{\scrR}(N,X)\in\CM\scrR\}.
\]
In the special case when $\scrR$ is isolated, the right hand side is equivalent to the condition
\[
\add N=\{ X\in\CM \scrR \mid \Ext^1_{\scrR}(N,X)=0\},
\]
which is, by definition, the notion of $N$ being a \emph{cluster tilting} object in the category $\uCM\scrR$.  This is the key link to cluster theory.  However, whilst helpful, as outlined in the introduction this analogy only takes us so far: $\scrR$ need not be isolated, and need not admit an NCCR (equivalently, a crepant resolution).  The general case is covered by using the language of MM modules.

\subsection{ Auslander--McKay correspondence}\label{AM section}
Let $\scrR$ be an arbitrary cDV singularity.  To any fixed minimal model $\scrX\to\Spec\scrR$, via slicing \S\ref{Pinkham sect} associate $\scrJ\subseteq\Updelta$, and then via \S\ref{int arr section} associate a finite arrangement $\scrH$ and an infinite $\scrH^{\aff}$.  These turn out to be independent of the choice of the minimal model, albeit in a slightly subtle way (see \cite[4.2.2]{IW9}).

As notation, for $N\in\CM\scrR$, set $\uEnd_\scrR(N)=\End_\scrR(N)/[\scrR]$, where $[\scrR]$ is the two-sided ideal of morphisms $N\to N$ which factor through $\add\scrR$.
\begin{thm}\label{AM main}
Let $\scrR$ be \textnormal{cDV}, then there exists is a bijection
\[
\begin{array}{c}
\begin{tikzpicture}
\node (A) at (0.1,0) {$(\mathrm{MM}\,\scrR)\cap(\CM\scrR)$};
\node (B) at (6,0) {$
\{\mbox{minimal models $f_i\colon\scrX_i\to\Spec \scrR$} \}
$};
\draw[<->] (1.85,0) -- node [above] {$\scriptstyle $} (3,0);
\node at (0,0.175) {$\phantom -$};
\end{tikzpicture}
\end{array}
\]
where the left-hand side is taken up to isomorphism, and the right-hand side is taken up to isomorphism of the $\scrX_i$ compatible with the morphisms $f_i$.  Under this bijection:
\begin{enumerate}
\item\label{AM main 1} For any \textnormal{MM} generator, its non-free indecomposable summands are in one-to-one correspondence with the exceptional curves in the corresponding minimal model.
\item For any fixed  \textnormal{MM} generator $N$, the quiver of $\uEnd_\scrR(N)$ encodes the dual graph of the corresponding minimal model.
\item\label{AM 3} The mutation graph of the  \textnormal{MM} generators coincides with the flops graph of the minimal models.  
\item The graphs in \eqref{AM 3} can be realised, inside the skeleton of $\scrH$, as a connected path. In this path every minimal model of $\scrR$ appears precisely once, and furthermore each wall crossing in this path corresponds to a flop.  When $\scrR$ is isolated, this path covers the whole of $\scrH$.
\end{enumerate}
\end{thm}

We remark that the bijective maps in \ref{AM main} are in fact explicit. The passage from left to right takes a given $N$ and associates a certain moduli space of representations of $\End_{\scrR}(N)$.  The passage from right to left takes global sections of the Van den Bergh tilting generator of $0$-perverse sheaves; we refer the reader to \cite[\S4.2]{HomMMP} for full details.

\medskip
In many ways, the strength in \ref{AM main} comes from the flexibility in being able to enlarge the left hand side of the bijection, and instead consider the set  $\mathrm{MM}\,\scrR$.  This moves us beyond what is seen by the birational geometry, even when $\scrR$ is isolated.  The following is in many ways the main result of \cite{IW9}.  As above, the infinite hyperplane  arrangement $\scrH^{\aff}$ associated to $\scrR$ strictly speaking depends on a choice of minimal model (or choice of maximal modifying module), but it turns out that this choice does not matter.

\begin{thm}
If $\scrR$ is \textnormal{cDV}, then there is an injection
\begin{equation}
\mathrm{MM}\,\scrR\hookrightarrow\{ \textnormal{chambers of }\scrH^{\aff}\}\label{inj aff}
\end{equation}
and furthermore the following statements hold.
\begin{enumerate}
\item If $N,L\in\mathrm{MM}\,\scrR$ with $L\ncong N$, then $L$ and $N$ are connected by a mutation at an indecomposable summand if and only if their chambers are adjacent.
\item When $\scrR$ is isolated, the map \eqref{inj aff} is a bijection.
\end{enumerate}
\end{thm}

In particular, this gives a \emph{complete} classification of noncommutative resolutions, and their variants such as MMAs, for cDV singularities. We remark (see \cite[\S9]{IW9}) that there are variations of the above results for arbitrary modifying modules, and this generality is required for applications to autoequivalences in \S\ref{autostab section}.  

For now, the following corollary summarises the above in the case when $\scrR$ is isolated. From a cluster perspective, the main insight is to move away from finite type by replacing the Frobenius category $\CM\scrR$ with the larger category $\refl\scrR$. Whilst this brings immediate benefits, allowing us to see the `affine' version of the theory, the category $\refl\scrR$ does not have an obvious cluster interpretation. 

\begin{cor}\label{bijection main isol}
If $\scrR$ is isolated \textnormal{cDV}, then there are bijections
\[
\begin{array}{rcl}
\textnormal{maximal rigid objects in }\CM\scrR&\longleftrightarrow&\{ \textnormal{chambers of }\scrH\},\\
\textnormal{maximal modifying objects in }\refl\scrR&\longleftrightarrow&\{ \textnormal{chambers of }\scrH^{\aff}\}.\\
\end{array}
\]
In both cases, wall crossing corresponds to mutation. 
\end{cor}

\begin{example}\label{ex: Haff 1}
To illustrate the second bijection in \ref{bijection main isol}, the left hand side below  illustrates an example of an $\scrH^{\aff}$ for a particular $cE_8$ singularity, with a dot drawn in each chamber.  The right hand side illustrates the mutation graph, where we connect two dots if and only if their chambers are adjacent.
\[
\includegraphics[angle=0,scale = 0.4]{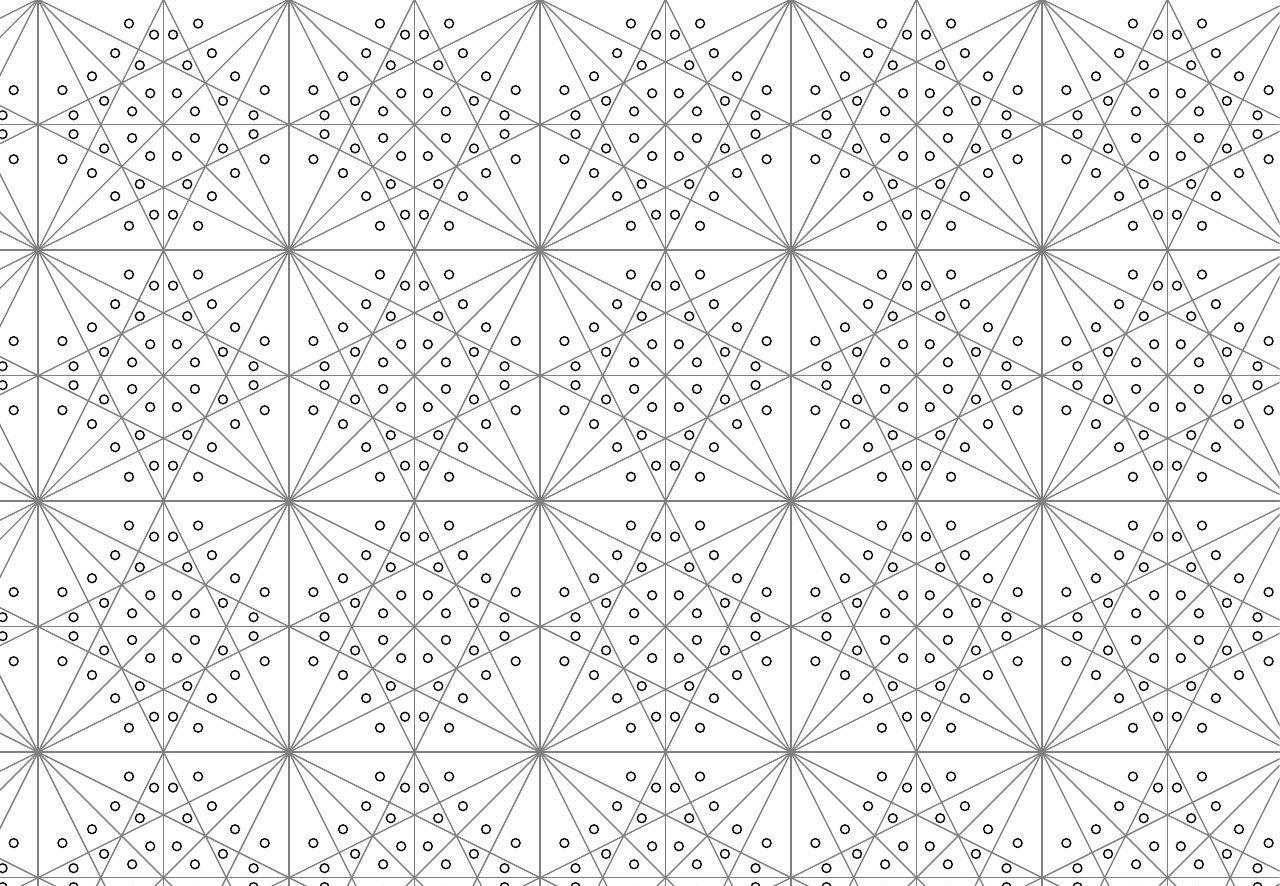}
\quad
\includegraphics[angle=0,scale = 0.4]{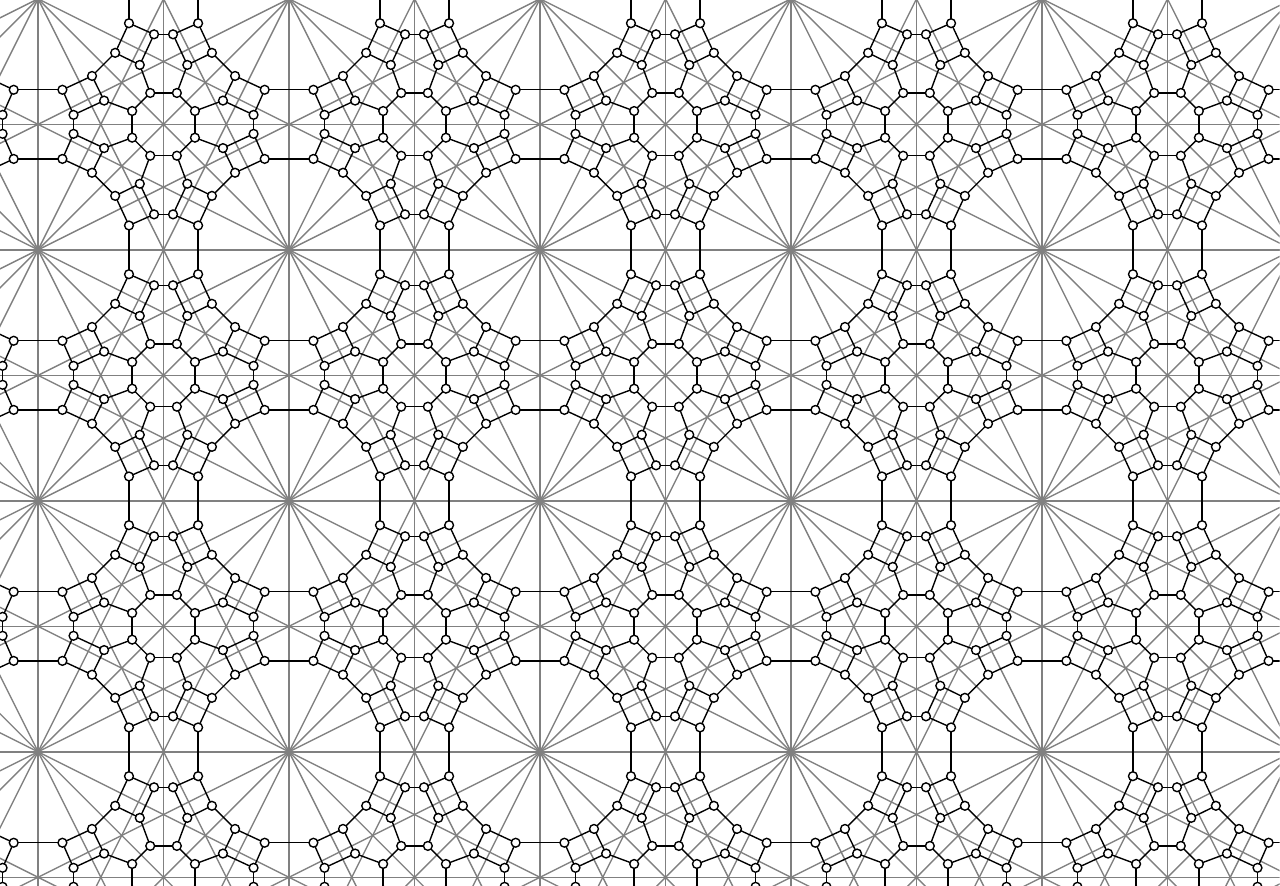}
\]
\end{example}

\section{Contraction Algebras}

Given a crepant partial resolution $\scrX\to\Spec\scrR$ with $\scrR$ cDV, this section explains how to construct the contraction algebra $\CA$.  When $\scrR$ is isolated this algebra is finite dimensional over $\mathbb{C}$ (i.e.\ GK-dimension zero), else it has GK-dimension one.  The construction is in fact much more general, but for ease throughout we mostly restrict to the cDV case.

\subsection{Deformation Theory Overview}\label{def theory section}
Contraction algebras are defined first and foremost using deformation theory. Let $X$ be a variety, and $\scrF\in\Qcoh X$.  To deform $\scrF$ via the classical Grothendieck formulation, consider the functor 
\[
\cDef\colon \cart_1\to\Sets
\]
from commutative local artinian $\mathbb{C}$-algebras to sets, sending
\[
(R,\m)\mapsto \left. \left \{ 
(\scrG,\mathsf{f})
\left|\begin{array}{l}
\scrG\in\Qcoh (X\times \Spec R),\\ 
\scrG\mbox{ is flat over }\Spec R,\\ 
\mathsf{f}\colon \scrG|_{X\times 0}\xrightarrow{\sim}\scrF.
\end{array}\right. \right\} \middle/ \sim \right.
\]
where $\sim$ is a natural equivalence relation which can be ignored for now.  The point is that, roughly speaking, a deformation of $\scrF$ is a sheaf $\scrG$ which is infinitesimally thicker than $\scrF$ (namely, a sheaf on $X\times \Spec R$), roughly of the same size (flatness), which recovers the original $\scrF$ (restricted to the zero fibre).

The above will referred to as the commutative deformation functor. Under nice, general situations, the functor $\cDef$ is \emph{prorepresentable}. For our purposes later (see \ref{contraction theorem}), this will not be enough, and we will require more general versions.

Let $\scrA$ be either $\Mod\Lambda$ or $\Qcoh Y$, where $\Lambda$ is a $\mathbb{C}$-algebra and $Y$ is a quasi-projective $\mathbb{C}$-scheme.  Pick $\scrF\in\scrA$, then the noncommutative deformation functor is defined to be
\[
\Def\colon \art_1\to\Sets
\]
from finite dimensional augmented $\mathbb{C}$-algebras to sets, sending
\[
(\Gamma,\n)\mapsto \left. \left \{ 
(\scrG,\upvartheta,f)
\left|\begin{array}{l}
\scrG\in\scrA,\\
\upvartheta\colon \Gamma\to\End_{\scrA}(\scrG)\mbox{ homomorphism},\\ 
-\otimes_\Gamma \scrG\colon\mod\Gamma\to \scrA \mbox{ is exact},\\ 
f\colon (\Gamma/\n)\otimes_\Gamma\scrG\xrightarrow{\sim}\scrF.
\end{array}\right. \right\} \middle/ \sim \right.
\]
where, for any given $\Gamma$, $(\scrG,\upvartheta,f)\sim (\scrG',\upvartheta',f')$ if and only if there exists $\uptau \colon \scrG\xrightarrow{\sim}\scrG'$ such that for all $r\in\Gamma$
\[
\begin{array}{c}
\begin{tikzpicture}
\node (a1) at (0,0) {$\scrG$};
\node (a2) at (1.5,0) {$\scrG'$};
\node (b1) at (0,-1.5) {$\scrG$};
\node (b2) at (1.5,-1.5) {$\scrG'$};
\draw[->] (a1) -- node[above] {$\scriptstyle \uptau$} (a2);
\draw[->] (a1) -- node[left] {$\scriptstyle \upvartheta(r)$} (b1);
\draw[->] (a2) -- node[right] {$\scriptstyle \upvartheta'(r)$} (b2);
\draw[->] (b1) -- node[above] {$\scriptstyle \uptau$} (b2);
\end{tikzpicture}
\end{array}
\quad
\mbox{and}
\quad
\begin{array}{c}
\begin{tikzpicture}
\node (a1) at (0,0) {$(\Gamma/\n)\otimes_\Gamma \scrG$};
\node (a2) at (3,0) {$(\Gamma/\n)\otimes_\Gamma \scrG'$};
\node (b) at (1.5,-1.5) {$\scrF$};
\draw[->] (a1) -- node[above] {$\scriptstyle 1\otimes \kern 1pt \uptau$} (a2);
\draw[->] (a1) -- node[pos=0.6,left,inner sep=8,anchor=east] {$\scriptstyle f$} (b);
\draw[->] (a2) -- node[pos=0.6,right,inner sep=8,anchor=west] {$\scriptstyle f'$} (b);
\end{tikzpicture}
\end{array}
\]
commute.  In the case $\scrA=\Qcoh X$, the commutative deformation functor is recovered as 
\[
\cDef=\Def|_{\cart_1}.
\]
The fact that the noncommutative deformation functor has \emph{more} test objects means that any universal sheaf for $\Def$ can only be \emph{larger} than, or equal to, the universal sheaf for $\cDef$.  In the situation of curves, this contains strictly more information.

\subsection{Deforming multiple curves}\label{NC multiple sction}
Noncommutative deformations give strictly more information than commutative deformations, even in the case of a single irreducible flopping curve.  However, it is in the situation of multiple intersecting curves where noncommutative deformations really come into their own. Indeed, it is precisely the Ext information \emph{between} curves that records their intersection data.  Throwing this away or merging it, which is what commutative deformation theory does, is at best counterproductive.

There are various equivalent ways of noncommutatively deforming curves. Here we approach the problem using DG-algebras. The main benefit of this approach is its brevity, and the main cost is the loss of all conceptual understanding. The case of a single curve, avoiding DG-algebras, is given in \cite[\S2]{DW1}, using the formulation in \S\ref{def theory section}.  More generally, multi-pointed deformations avoiding DG-algebras are explained in \cite{Laudal, Eriksen, Eriksen2,Kawamata,Kawamata2}.

Composing arrows as in the conventions, recall that a DG~category is a graded category $\mathsf{A}$ whose morphism spaces are endowed with a differential~$\updelta$, i.e.\ homogeneous maps of degree one satisfying $\updelta^2=0$, such that 
\[
\updelta(gf)=g(\updelta f)+(-1)^t(\updelta g)f
\]
for all $g\in\Hom_{\mathsf{A}}(a,b)$ and all $f\in\Hom_{\mathsf{A}}(b,c)_t$ for $t \in \mathbb{Z}$.  For our purposes here, it will be important to study the category $\DG_n$, which has as objects those DG categories with precisely $n$ objects. Similarly, write $\art_n$ for the category of $n$-pointed artinian $\mathbb{C}$-algebras.

\begin{notation}
Suppose that $(\mathsf{A},\updelta)\in\DG_n$, and $(\Gamma,\n)\in\art_n$, and set $\n_{ij}:=e_i\n\hspace{0.1em}e_j$.  Define $\mathsf{A}\uotimes\n:=\bigoplus_{i,j=1}^n(\mathsf{A}\uotimes\n)_{ij}$, where
\[
(\mathsf{A}\uotimes\n)_{ij}:=\Hom_{\mathsf{A}}(i,j)\otimes_{\K} \n_{ij}.
\]
Observe that $\mathsf{A}\uotimes\n$ has the natural structure of an object in $\DG_n$ (but with no units) with differential $d(a\otimes x):=d(a)\otimes x$. Thus we may consider $\mathsf{A}\uotimes\n$ as a DGLA, with bracket
\[
\quad[a\otimes x,b\otimes y]:=ab\otimes xy-(-1)^{\deg(a)\deg(b)}ba\otimes yx
\] 
for homogeneous $a,b\in\mathsf{A}$.
\end{notation}

\begin{defin}
Given $(\mathsf{A},\updelta)\in\DG_n$,  the associated DG deformation functor
\[
\Def^{\mathsf{A}}\colon\art_n\to\Sets
\]
is defined by sending
\[
(\Gamma,\n)\mapsto
\left. \left \{ 
\upxi\in \mathsf{A}^1\underline{\otimes}\,\n
\left|
\,
\updelta(\upxi)+\tfrac{1}{2}[\upxi,\upxi]=0
\right. \right\} \middle/ \sim \right.
\]
where as usual the equivalence relation $\sim$ is induced by the gauge action.  Explicitly, two elements $\upxi_1,\upxi_2\in \mathsf{A}^1\uotimes\n$ are said to be gauge equivalent if there exists $x\in\mathsf{A}^0\uotimes\n$ such that
\[
\upxi_2=e^x*\upxi_1:=\upxi_1+\sum_{j=0}^{\infty}\frac{([x,-])^j}{(j+1)!}([x,\upxi_1]-d(x)).
\]
\end{defin}

Under very general situations, the functor $\Def^{\mathsf{A}}$ is \emph{prorepresentable}.  We shall not use this, as the prorepresenting object in this construction is difficult to control.  In our cDV setting $\scrX\to\Spec\scrR$ below,  we will instead build the prorepresenting object directly.

\medskip
Given a crepant partial resolution $f\colon\scrX\to\Spec\scrR$ with $\scrR$ cDV, consider the curves $\bigcup_{i=1}^n\mathrm{C}_i$ above the origin, with each $\mathrm{C}_i\cong\mathbb{P}^1$, and the resulting structure sheaves $\{ \scrO_{\mathrm{C}_i}\}_{i=1}^n$.  For each $i$, choose an injective resolution $\scrO_{\mathrm{C}_i}(-1)\hookrightarrow\scrI_i$, and form the endomorphism DG-algebra $\sEnd(\bigoplus_{i=1}^n\scrI_i)$.  Crucially, we can view this in $\DG_n$, and so applying the above generalities obtain a noncommutative deformation functor, written
\[
\Def_{f}\colon\art_n\to\Sets.
\]
On the other hand, exactly as in \ref{AM main}, for $f\colon\scrX\to\Spec\scrR$ write $N$ for the global sections of the Van den Bergh tilting generator of $0$-perverse sheaves, and  set $\CA=\uEnd_\scrR(N)$.  The following then relates deformation theory to cluster theory.
\begin{thm}\label{prorep main}
With notation as above, $\Def_{f}\cong\Hom(\CA,-)$.
\end{thm}
We call $\CA=\uEnd_\scrR(N)$ the \emph{contraction algebra} associated to $f\colon\scrX\to\Spec\scrR$.  The pro-representing object is always unique, up to isomorphism, and so in fact any algebra in this isomorphism class will do.  We remark  that in addition to deformation theory and cluster theory, there is a third method to obtain the contraction algebra, as the factor of an algebra obtained via tilting.  This third viewpoint becomes important in \S\ref{autostab section}.

\medskip
It is both easy and desirable to tweak the above to work for any \emph{choice} of subset of curves.  Indeed, for any subset $\bigcup_{i\in I}\mathrm{C}_i$ of curves, there is a corresponding deformation functor $\Def_{I}$, built using $\sEnd(\bigoplus_{i\in I}\scrI_i)\in\DG_{|I|}$.  It turns out that in this more general setting,
\begin{equation}
\Def_{I}\cong\Hom(\End_{\scrR}(N)/\scrK_I,-),\label{cont general}
\end{equation}
where $\scrK_I$ is the two-sided ideal consisting of morphisms $N\to N$ which factor through a certain summand of $N$; see \cite{DW3} for full details.  The larger the choice of curves, the smaller the ideal $\scrK_I$, and thus the larger the algebra $\End_{\scrR}(N)/\scrK_I$.  It is an easy fact that all the $\End_{\scrR}(N)/\scrK_I$ are factors, by an idempotent, of the contraction algebra $\CA$ associated to the full contraction $f$. In particular, each is \emph{very} easy to calculate, provided that $\CA$ is known explicitly.

Given a choice of curves $\bigcup_{i\in I}\mathrm{C}_i$, as in \S\ref{min models sect} contract them to yield morphisms
\[
\scrX\xrightarrow{g_I}\scrX_I\to\Spec\scrR.
\]
One of the main points in HomMMP \cite{HomMMP} is that $\End_{\scrR}(N)/\scrK_I$ is the contraction algebra for the contraction $g_I$.  In particular, by \ref{contraction theorem} below the finite dimensionality of $\End_{\scrR}(N)/\scrK_I$ determines if $\bigcup_{i\in I}\mathrm{C}_i$ flops, or not.  Thus, determining which subsets of curves flop, and which do not, becomes much easier, as all this information is encoded in $\CA$ and its (easy) idempotent factors.

\begin{example}
Consider the $\scrR$ in \ref{7.7HomMMPa}, which is $cD_4$. As explained in \cite[7.7]{HomMMP}, the contraction algebra $\CA$ associated to one of the crepant resolutions $\scrX\to\Spec\scrR$ is given by the following quiver with relations.
\begin{equation}
\begin{array}{ccc}
\begin{array}{c}
\begin{tikzpicture}[bend angle=10, looseness=1,>=stealth]
\node (1) at (0,0) [vertex] {};
\node (1a) at ($(1)+(120:2pt)$)  {};
\node (1b) at ($(1)+(60:2pt)$)  {};
\node (3) at (2,0) [vertex] {};
\draw[->,black,bend left]  (3) to node[gap] {$\scriptstyle b$} (1);
\draw[->,black,bend left]  (1) to node[gap] {$\scriptstyle B$} (3);
\draw[->]  (1) edge [in=100,out=160,loop,looseness=10] node[left]{$\scriptstyle \ell$} (1);
\draw[->]  (1) edge [in=260,out=200,loop,looseness=10] node[left]{$\scriptstyle z$} (1);
\end{tikzpicture}
\end{array}
&&
{\scriptsize
\begin{array}{l}
Bb=2\ell^2\\
bz=0\\
zB=0\\
\ell z=-z\ell
\end{array}
}
\end{array}\label{full cont}
\end{equation}
There are two curves above the origin, which correspond to the vertices. It is clear that $\CA$ is infinite dimensional.  Choosing the curve which corresponds to the left hand vertex, we now test whether it flops.  To do this, we delete the \emph{other} vertex, and all arrows that pass through it, leaving
\[
\begin{array}{ccc}
\begin{array}{c}
\begin{tikzpicture}[bend angle=10, looseness=1,>=stealth]
\node (1) at (0,0) [vertex] {};
\node (1a) at ($(1)+(120:2pt)$)  {};
\node (1b) at ($(1)+(60:2pt)$)  {};
\node (3) at (2,0) [vertexblank] {};
\draw[black!40!white,->,bend left]  (3) to node[gap] {$\scriptstyle b$} (1);
\draw[black!40!white,->,bend left]  (1) to node[gap] {$\scriptstyle B$} (3);
\draw[->]  (1) edge [in=100,out=160,loop,looseness=10] node[left]{$\scriptstyle \ell$} (1);
\draw[->]  (1) edge [in=260,out=200,loop,looseness=10] node[left]{$\scriptstyle z$} (1);
\end{tikzpicture}
\end{array}
&&
{\scriptsize
\begin{array}{l}
0=2\ell^2\\
\ell z=-z\ell.
\end{array}
}
\end{array}
\]
This is the contraction algebra associated to the contraction of the left hand curve, which is evidently infinite dimensional, and thus the left hand curve does not flop.  To test whether the right hand curve flops, instead delete the left hand vertex in \eqref{full cont} and all arrows passing through it, leaving
\[
\begin{array}{ccc}
\begin{array}{c}
\begin{tikzpicture}[bend angle=10, looseness=1,>=stealth]
\node (1) at (0,0) [vertexblank] {};
\node (1a) at ($(1)+(120:2pt)$)  {};
\node (1b) at ($(1)+(60:2pt)$)  {};
\node (3) at (2,0) [vertex] {};
\draw[black!40!white,->,bend left]  (3) to node[gap] {$\scriptstyle b$} (1);
\draw[black!40!white,->,bend left]  (1) to node[gap] {$\scriptstyle B$} (3);
\draw[black!40!white,->]  (1) edge [in=100,out=160,loop,looseness=10] node[left]{$\scriptstyle \ell$} (1);
\draw[black!40!white,->]  (1) edge [in=260,out=200,loop,looseness=10] node[left]{$\scriptstyle z$} (1);
\end{tikzpicture}
\end{array}
&&
{\scriptsize
\begin{array}{l}
\phantom{0=2\ell^2}\\
\phantom{\ell z=-z\ell}
\end{array}
}
\end{array}
\]
This is just the complex numbers, which is certainly finite dimensional.  In fact, this right hand curve gives the Atiyah flop.
\end{example}

\subsection{Main Results and Conjectures}\label{main contalg subsection}
The contraction algebra $\CA$ admits three compatible structures.  The first, the structure of an algebra, conjecturally gives classification.  The second is the structure of an $\scrR$-module, whose support gives rise to the Contraction Theorem below. The third is the structure of an $\End_\scrR(\scrR\oplus N)$-module, which via \eqref{VdBtilt} gives a sheaf $\scrE\in\coh\scrX$.  Stalk locally, this is the universal sheaf from deformation theory, up to additive equivalence.  From this third viewpoint, the key point is that $\scrE$ is a perfect complex when $\scrX$ is only mildly singular (see \ref{perf and triangle}), which should be viewed as the `hidden smoothness' needed for the derived equivalence technology to work well.

It is the ability to view $\CA$ as an $\scrR$-module that allows the contracted locus to be determined. The following has a global version in \cite{DW2}, and is true much more generally, but here we again restrict to the local cDV setting for simplicity.

\begin{thm}[Contraction Theorem, \cite{DW2}]\label{contraction theorem}
Suppose that $f\colon \scrX\to\Spec\scrR$ is a crepant partial resolution, and write $Z$ for the locus in $\Spec\scrR$ for which $f$ is not an isomorphism.  Then $\Supp_{\scrR}\CA=Z$.  In particular, the following are equivalent
\begin{enumerate}
\item $f$ is a flopping contraction (namely, we are in situation \ding{192}, not \ding{193} in \eqref{two cases intro}) 
\item $\dim_{\mathbb{C}}\CA<\infty$.
\end{enumerate}
\end{thm}

In contrast, it is the structure of an algebra on $\CA$ which conjecturally classifies.

\begin{conj}[The Classification Conjecture \cite{DW1}]\label{class:conj} 
Let $\scrX_1\to\Spec \scrR_1$ and $\scrX_2\to\Spec \scrR_2$ be $3$-fold flopping contractions, where $\scrX_i$ are smooth, and $\scrR_i$ are complete local.  Denote their corresponding contraction algebras by $\Acon$ and $\Bcon$, respectively.  Then
\[
\scrR_1\cong \scrR_2\iff\Db(\mod\Acon)\simeq\Db(\mod\Bcon).
\]
\end{conj}
The derived equivalence condition on the right hand side can be greatly simplified, and the classification boils down to isomorphism classes of algebras.  Indeed, by work of August \cite[1.5]{August1}, the right hand side is equivalent to $\Acon$ being \emph{isomorphic} to an iterated mutation of $\Bcon$, of which there are only finitely many.  Thus, the Classification Conjecture is really a statement about isomorphism classes of contraction algebras, and not about derived categories.  In the case of irreducible flops, $\Acon$ is local, and so there are no further basic algebras in its derived equivalence class, up to isomorphism. Thus in this case we rid ourselves of iterated mutations, and the conjecture reduces to the statement that irreducible smooth $3$-fold flopping contractions are classified by their contraction algebras, up to isomorphism.

There is now a body of evidence underpinning \ref{class:conj}.  Type $A$ was known from the original work of Reid \cite{Pagoda}, and Type $D$ evidence is provided by \cite{BW, Okke, Kawamata3}.  Various enhanced versions of the conjecture are known to be true \cite{Hua, HuaKeller, Booth}, but indeed the point of the conjecture is that these enhanced structures are not necessary. The contraction algebra remains the finest known curve invariant.

\medskip
Both from a future classification perspective, and for purely algebraic reasons, it is important to know which algebras can be realised as (more precisely, are isomorphic to) the contraction algebra of some crepant resolution $\scrX\to\Spec\scrR$, with $\scrR$ cDV.  Below, such algebras are called \emph{geometric}.  As notation to state the following, for $f\in\mathbb{C}\llangle x,y\rrangle$ consider the associated Jacobi algebra
\[
\Jac(f)\colonequals \frac{\mathbb{C}\llangle x,y\rrangle}{\llbr\updelta_x f,\updelta_y f\rrbr}
\]
where $\updelta_x, \updelta_y$ are the cyclic derivatives, and $\llbr I\rrbr$ denotes the closure of the two sided ideal $I$. The following covers the case when there is only one curve above the origin: more general statements exist, but are harder to state.

\begin{conj}[Realisation Problem, \cite{BW2}] \label{realisation:conj}
If $f\in\mathbb{C}\llangle x,y\rrangle$ satisfies $\GKdim\Jac(f)\leq 1$, then $\Jac(f)$ is geometric. In particular every finite dimensional superpotential algebra
\[
\frac{\mathbb{C}\llangle x,y\rrangle}{\llbr\updelta_x f,\updelta_y f\rrbr}
\]
can be constructed as the contraction algebra of some irreducible $3$-fold flop.
\end{conj}
Of the two conjectures,  \ref{realisation:conj} is by far the most outrageous. It arose not from any great philosophical underpinning, but from original disbelief that was superseded, over a number of years, by significant evidence from extensive computer algebra searches, theoretical noncommutative advances \cite{BW2, Iyudu, IS, IS2}, and from DT theory \cite{Davison}.  The paper \cite{BW2} contains a summary of the best known results in this direction.

\subsection{Curve Counting}\label{Toda section}
The contraction algebra $\CA$ contains all known curve invariants, albeit extracting them explicitly can be challenging.  Here, we restrict to smooth irreducible flopping contractions $f\colon\scrX\to\Spec\scrR$, which have an associated tuple of integers $(n_1,\hdots,n_\ell)$ called the \emph{Gopakumar--Vafa} invariants,  defined below.  The idea is to deform $f$ into a disjoint union of the simplest type of flopping curves, and to count the curves there. The subtlety comes from the fact that we don't just naively count the total number of such curves in the deformation, instead we use the information of the flat family to split (or, refine) the count into the number of curves with a given curve class. 

More precisely, by \cite[\S2.1]{BKL} there exists a flat deformation
\[
\begin{tikzpicture}[>=stealth]
\node (X) at (0,0) {$\cX$};
\node (Y) at (0,-1) {$\cY$};
\node (T) at (0.75,-1.75) {$T$};
\draw[->] (X)--(Y);
\draw[->] (Y)--(T);
\draw[->,densely dotted] (X)--(T);
\end{tikzpicture}
\]
for some Zariski open neighbourhood $T$ of $0\in\mathbb{A}^1$, such that the central fibre $g_0\colon X_0\to Y_0$ is isomorphic to $f\colon\scrX\to\Spec\scrR$, and further all other fibres $g_t\colon X_t\to Y_t$ for $t\in T\backslash \{0\}$  are flopping contractions whose exceptional locus is a disjoint union of $(-1,-1)$-curves.

Regarding the flopping curve $\mathrm{C}$ of $f\colon\scrX\to\Spec\scrR$ as a curve in $\mathcal{X}$, then the GV invariant $n_j$ is defined to be the number of $g_t$-exceptional $(-1,-1)$-curves $\mathrm{C}'$ with curve class $j[\mathrm{C}]$, namely for every line bundle $\scrL$ on $\cX$, 
\[
\deg (\scrL|_{\mathrm{C}'})=j \deg(\scrL|_\mathrm{C}).
\] 
Toda's formula relates these numbers to the dimension of the contraction algebra.  As notation, write $\CA^{\ab}$ for the abelianisation of $\CA$, which is the algebra obtained from $\CA$ by further factoring by the ideal of all possible commutators.  Then for $f\colon \scrX\to \Spec \scrR$ a smooth length $\ell$ irreducible flopping contraction, Toda \cite{TodaGV} shows that $n_1=\dim_{\mathbb{C}}\CA^{\ab}$ and furthermore
 \[
 \dim_{\mathbb{C}}\CA=\dim_{\mathbb{C}}\CA^{\ab}+\sum_{j=2}^{\ell}j^2\cdot n_j.
 \]
There is a more general formula that gives the dimension of the contraction algebra in terms of the GV invariants for arbitrary (multi-curve) smooth flopping contractions \cite{TodaUtah}, but this is mildly more technical to state.  Further, Hua--Toda \cite[4.6]{HuaToda} show that the GV invariants are a property of the isomorphism class of the contraction algebra in the case of smooth irreducible flops (which is to be expected from \ref{class:conj}), however it is still open how to easily extract them.  

The numerical GV invariants are precisely that, they are numerical information.  The algebra structure on the contraction algebra is \emph{additional} information:  there are examples of non-isomorphic flops that admit the same GV invariants \cite{BW}.  In this case, it is the algebra structure on the contraction algebra that distinguishes them.

\section{Extended Example: Type A}\label{Type A section}

In this section we consider any choice of  $f\in\mathbb{C}\llsq x,y\rrsq$, and set 
\[
\scrR\colonequals \mathbb{C}\llsq u,v,x,y\rrsq/(uv-f).
\]   
These all turn out to be $cA_m$, where $m=\mathrm{ord}(f)-1$.  Such singularities have been studied since Viehweg \cite{Viehweg}, with the homological and noncommutative properties being developed in increasing levels of generality in \cite{Nag,BIKR,DH,IW5,IW6}.  

In this section we review the state of the art, and give a combinatorial model for flopping curves in the minimal models of $\Spec\scrR$.  

\subsection{Maximal Rigid Objects} 
With $\scrR$ as above, the fact that $\scrR$ is cDV, in fact $cA_m$ with $m=\mathrm{ord}(f)-1$ is somewhat implicit in the MMP literature; a full proof is given in \cite[6.1(e)]{BIKR}.  What makes understanding the modifying modules of $\scrR$ possible in this context is that \emph{all} modifying modules are a direct sum of rank one reflexive modules. This is a direct consequence of $\scrR$ slicing to a type $A$ Kleinian singularity, and is far from typical behaviour.

Luck being what it is, all modifying modules can be described by remaining entirely within the realm of the class group $\Cl(\scrR)$, and taking direct sums of elements thereof.  The following theorem was established for those isolated $\scrR$ which furthermore admit an NCCR in \cite{BIKR}, then again by purely algebraic means in \cite{DH}.  The isolated assumption, and the assumption on admitting an NCCR, were both stripped in \cite{IW5,IW6}.  

As notation, write $f=f_1\hdots f_n$ with each $f_i$ irreducible, and for any element of the symmetric group $\upomega\in\mathfrak{S}_n$ consider the following $\scrR$-module
\[
M^\upomega=\scrR\oplus (u,f_{\upomega(1)})\oplus (u,f_{\upomega(1)}f_{\upomega(2)})\oplus\hdots\oplus (u,f_{\upomega(1)}\hdots f_{\upomega(n-1)}),
\]
where $(u,g)$ denotes the ideal of $\scrR$ generated by $u$ and $g$.  If there are repetition in the $f_i$, note that the $M^{\upomega}$ are not all distinct.

\begin{thm}\label{IWclassification}
Let $\scrR\colonequals \mathbb{C}\llsq u,v,x,y\rrsq/(uv-f)$, then the following hold.
\begin{enumerate}
\item\label{IWclassification 1} The basic elements of  $(\mathrm{MM}\,\scrR)\cap(\CM\scrR)$ are precisely the $M^\upomega$, with $\upomega\in\mathfrak{S}_n$.  
\item\label{IWclassification 2} The basic \textnormal{MM} $\scrR$-modules are precisely $(I\otimes_R M^\upomega)^{**}$ with $\upomega\in\mathfrak{S}_n$, and $I\in\Cl(R)$.
\end{enumerate}
\end{thm}

Grouping together irreducible terms and stripping units, we may write $f=g_1^{a_1}\hdots g_{t}^{a_t}$ where each $g_i$ is irreducible and $(g_i)\neq(g_j)$ for $i\neq j$.  In this notation it is easy to see from \ref{IWclassification}\eqref{IWclassification 1} that there are precisely $\frac{(a_1+\hdots + a_t)!}{a_1!\hdots a_t!}$ elements of $(\mathrm{MM}\,\scrR)\cap(\CM\scrR)$. By \ref{AM main} there are thus that number of minimal models, however this is overkill: in this case, it is possible to deduce this much more simply, by using the combinatorial models below.

\subsection{Minimal Models}\label{min model A sect}
Producing minimal models of $\scrR\colonequals \mathbb{C}\llsq u,v,x,y\rrsq/(uv-f)$ using elements $\upomega\in\mathfrak{S}_n$ is also well known and implicit in the literature. Indeed, for such an $\upomega$, blowing up $\Spec\scrR$ at the ideal $(u,f_{\upomega(1)})$ yields a scheme whose fibre above the unique closed point is a single $\mathbb{P}^1$ with two singularities, which suitably locally are $uv=f_{\upomega(1)}$, and $uv=\frac{f}{f_{\upomega(1)}}$. We sketch this as follows.
\begin{center}
\begin{tikzpicture}[xscale=0.6,yscale=0.6]
\draw[black] (-0.1,-0.04,0) to [bend left=25] (2.1,-0.04,0);
\filldraw [red] (0,0,0) circle (1pt);
\filldraw [red] (2,0,0) circle (1pt);
\node at (0,-0.4,0) {$\scriptstyle f_{\upomega(1)}$};
\node at (2,-0.4,0) {$\scriptstyle f/f_{\upomega(1)}$};
\end{tikzpicture}
\end{center}
 Proceeding by induction, we iteratively blowup $(u,f_{\upomega(2)})$ in the second chart, etc, and after finally many steps obtain $\scrX^\upomega\to\Spec\scrR$, where the fibre above the origin is a chain of curves
 \begin{center}
\begin{tikzpicture}[xscale=0.7,yscale=0.6]
\draw[black] (-0.1,-0.04) to [bend left=25] (2.1,-0.04,0);
\draw[black] (1.9,-0.04) to [bend left=25] (4.1,-0.04,0);
\draw[black] (6.9,-0.04) to [bend left=25] (9.1,-0.04,0);
\draw[black] (8.9,-0.04) to [bend left=25] (11.1,-0.04,0);
\filldraw [red] (0,0) circle (1pt);
\filldraw [red] (2,0) circle (1pt);
\filldraw [red] (4,0) circle (1pt);
\filldraw [red] (7,0) circle (1pt);
\filldraw [red] (9,0) circle (1pt);
\filldraw [red] (11,0) circle (1pt);
\node at (0,-0.4) {$\scriptstyle f_{\upomega(1)}$};
\node at (2,-0.4) {$\scriptstyle f_{\upomega(2)}$};
\node at (4,-0.4) {$\scriptstyle f_{\upomega(3)}$};
\node at (5.5,0) {$\scriptstyle \hdots$};
\node at (7,-0.4) {$\scriptstyle f_{\upomega(n-2)}$};
\node at (9,-0.4) {$\scriptstyle f_{\upomega(n-1)}$};
\node at (11,-0.4) {$\scriptstyle f_{\upomega(n)}$};
\end{tikzpicture}
\end{center}
and where each red dot is suitably locally given by the equation $uv=f_{\upomega(i)}$.  Since each $f_i$ is irreducible, all these points are factorial, and so $\scrX^\upomega\to\Spec\scrR$ is a minimal model. Full details of this construction are given in \cite[\S5]{IW5}.

\begin{example}
Consider again the suspended pinch point $\scrR=\mathbb{C}\llsq u,v,x,y\rrsq/(uv-x^2y)$ from \ref{toric ex}\eqref{toric ex 2}. Writing $x^2y=yxx=xyx=xxy$, then the minimal models are
\[
\begin{tikzpicture}[xscale=0.3,yscale=0.3]
\draw[black] (-0.1,-0.04) to [bend left=25] (2.1,-0.04,0);
\draw[black] (1.9,-0.04) to [bend left=25] (4.1,-0.04,0);
\filldraw [red] (0,0) circle (1pt);
\filldraw [red] (2,0) circle (1pt);
\filldraw [red] (4,0) circle (1pt);
\node at (0,-0.7) {$\scriptstyle y$};
\node at (2,-0.7) {$\scriptstyle x$};
\node at (4,-0.7) {$\scriptstyle x$};
\end{tikzpicture}
\qquad
\begin{tikzpicture}[xscale=0.3,yscale=0.3]
\draw[black] (-0.1,-0.04) to [bend left=25] (2.1,-0.04,0);
\draw[black] (1.9,-0.04) to [bend left=25] (4.1,-0.04,0);
\filldraw [red] (0,0) circle (1pt);
\filldraw [red] (2,0) circle (1pt);
\filldraw [red] (4,0) circle (1pt);
\node at (0,-0.7) {$\scriptstyle x$};
\node at (2,-0.7) {$\scriptstyle y$};
\node at (4,-0.7) {$\scriptstyle x$};
\end{tikzpicture}
\qquad
\begin{tikzpicture}[xscale=0.3,yscale=0.3]
\draw[black] (-0.1,-0.04) to [bend left=25] (2.1,-0.04,0);
\draw[black] (1.9,-0.04) to [bend left=25] (4.1,-0.04,0);
\filldraw [red] (0,0) circle (1pt);
\filldraw [red] (2,0) circle (1pt);
\filldraw [red] (4,0) circle (1pt);
\node at (0,-0.7) {$\scriptstyle x$};
\node at (2,-0.7) {$\scriptstyle x$};
\node at (4,-0.7) {$\scriptstyle y$};
\end{tikzpicture}
\]
These combinatorial models match the pictures, left to right, in \ref{toric ex}\eqref{toric ex 2}.  Since the points $uv=x$ and $uv=y$ are smooth, in this example all minimal models are smooth.
\end{example}

The link between the $\scrX^{\upomega}$ defined here and the $M^\upomega$ of the previous subsection is provided by the Van den Bergh tilting bundle, which gives a derived equivalence
\begin{equation}
\Db(\coh \scrX^{\upomega})\to
\Db(\mod\End_{\scrR}(M^{\upomega})).\label{VdB for Type A}
\end{equation}
This allows us to use derived equivalences between the varying $\End_{\scrR}(M^{\upomega})$, with $\upomega\in\mathfrak{S}_n$, to obtain equivalences between the various minimal models.

\subsection{Mutation Combinatorics}
Mutations exist for an arbitrary modifying module, but we restrict here to mutations between the $M^{\upomega}$, namely between the members of $(\mathrm{MM}\,\scrR)\cap(\CM\scrR)$, since these are in some sense the most birational. Even with this restriction, it is worth emphasising that whilst these mutations include \emph{all} the flop functors, in the non-isolated situation it also contains more.

Mutation is a process in which summand(s) of $M^{\upomega}$ are replaced by others, in a universal way.  Details and definitions can be found in \cite[\S6]{IW3}, but the setting here, the following combinatorial model suffices.  To simplify the discussion, consider the case $uv=f_1\hdots f_7$, and for $\upomega\in\mathfrak{S}_7$ write $g_i=f_{\upomega(i)}$.  Under the bijection in \ref{AM main}\eqref{AM main 1}, the summands of $M^\upomega$ correspond to the exceptional curves as illustrated below.
\begin{center}
\begin{tikzpicture}[xscale=0.9,yscale=0.6]
\draw[black] (-0.1,-0.04) to [bend left=25] (2.1,-0.04);
\draw[black] (1.9,-0.04) to [bend left=25] (4.1,-0.04);
\draw[black] (3.9,-0.04) to [bend left=25] (6.1,-0.04);
\draw[black] (5.9,-0.04) to [bend left=25] (8.1,-0.04);
\draw[black] (7.9,-0.04) to [bend left=25] (10.1,-0.04);
\draw[black] (9.9,-0.04) to [bend left=25] (12.1,-0.04);
\filldraw [red] (0,0) circle (1pt);
\filldraw [red] (2,0) circle (1pt);
\filldraw [red] (4,0) circle (1pt);
\filldraw [red] (6,0) circle (1pt);
\filldraw [red] (8,0) circle (1pt);
\filldraw [red] (10,0) circle (1pt);
\filldraw [red] (12,0) circle (1pt);
\node at (0,-0.4) {$\scriptstyle g_1$};
\node at (2,-0.4) {$\scriptstyle g_2$};
\node at (4,-0.4) {$\scriptstyle g_3$};
\node at (6,-0.4) {$\scriptstyle g_4$};
\node at (8,-0.4) {$\scriptstyle g_5$};
\node at (10,-0.4) {$\scriptstyle g_6$};
\node at (12,-0.4) {$\scriptstyle g_7$};
\node at (-1,1) {$\scriptstyle\scrR$};
\node at (0,1) {$\scriptstyle\oplus$};
\node at (1,1) {$\scriptstyle (u,g_1)$};
\node at (2,1) {$\scriptstyle\oplus$};
\node at (3,1) {$\scriptstyle (u,g_1g_2)$};
\node at (4,1) {$\scriptstyle\oplus$};
\node at (5,1) {$\scriptstyle (u,g_1\hdots g_3)$};
\node at (6,1) {$\scriptstyle\oplus$};
\node at (7,1) {$\scriptstyle (u,g_1\hdots g_4)$};
\node at (8,1) {$\scriptstyle\oplus$};
\node at (9,1) {$\scriptstyle (u,g_1\hdots g_5)$};
\node at (10,1) {$\scriptstyle\oplus$}; 
\node at (11,1) {$\scriptstyle (u,g_1\hdots g_6)$};
\end{tikzpicture}
\end{center}
Now, to mutate at the summand $(u,g_1g_2)\oplus (u,g_1\hdots g_3)\oplus (u, g_1\hdots g_5)$ say, first think of this as a subset $\scrJ$ of the curves and draw a dotted box around them.  Mutation then acts via the following intuitive picture:
\begin{equation}
\begin{array}{c}
\begin{tikzpicture}[xscale=0.9,yscale=0.6]
\draw[black] (-0.1,-0.04,0) to [bend left=25] (2.1,-0.04,0);
\draw[black] (1.9,-0.04,0) to [bend left=25] (4.1,-0.04,0);
\draw[black] (3.9,-0.04,0) to [bend left=25] (6.1,-0.04,0);
\draw[black] (5.9,-0.04,0) to [bend left=25] (8.1,-0.04,0);
\draw[black] (7.9,-0.04,0) to [bend left=25] (10.1,-0.04,0);
\draw[black] (9.9,-0.04,0) to [bend left=25] (12.1,-0.04,0);
\filldraw [red] (0,0,0) circle (1pt);
\filldraw [red] (2,0,0) circle (1pt);
\filldraw [red] (4,0,0) circle (1pt);
\filldraw [red] (6,0,0) circle (1pt);
\filldraw [red] (8,0,0) circle (1pt);
\filldraw [red] (10,0,0) circle (1pt);
\filldraw [red] (12,0,0) circle (1pt);
\node at (0,-0.4,0) {$\scriptstyle g_{1}$};
\node at (2,-0.4,0) {$\scriptstyle g_{2}$};
\node at (4,-0.4,0) {$\scriptstyle g_{3}$};
\node at (6,-0.4,0) {$\scriptstyle g_{4}$};
\node at (8,-0.4,0) {$\scriptstyle g_{5}$};
\node at (10,-0.4,0) {$\scriptstyle g_{6}$};
\node at (12,-0.4,0) {$\scriptstyle g_{7}$};
\draw [densely dotted] (2,-0.1,0) -- (6,-0.1,0) -- (6,0.5,0) -- (2,0.5,0) -- cycle;
\draw [densely dotted] (8,-0.1,0) -- (10,-0.1,0) -- (10,0.5,0) -- (8,0.5,0) -- cycle;
\draw[black] (-0.1,-3.04,0) to [bend left=25] (2.1,-3.04,0);
\draw[black] (1.9,-3.04,0) to [bend left=25] (4.1,-3.04,0);
\draw[black] (3.9,-3.04,0) to [bend left=25] (6.1,-3.04,0);
\draw[black] (5.9,-3.04,0) to [bend left=25] (8.1,-3.04,0);
\draw[black] (7.9,-3.04,0) to [bend left=25] (10.1,-3.04,0);
\draw[black] (9.9,-3.04,0) to [bend left=25] (12.1,-3.04,0);
\filldraw [red] (0,-3,0) circle (1pt);
\filldraw [red] (2,-3,0) circle (1pt);
\filldraw [red] (4,-3,0) circle (1pt);
\filldraw [red] (6,-3,0) circle (1pt);
\filldraw [red] (8,-3,0) circle (1pt);
\filldraw [red] (10,-3,0) circle (1pt);
\filldraw [red] (12,-3,0) circle (1pt);
\node at (0,-3.4,0) {$\scriptstyle g_{1}$};
\node at (2,-3.4,0) {$\scriptstyle g_{4}$};
\node at (4,-3.4,0) {$\scriptstyle g_{3}$};
\node at (6,-3.4,0) {$\scriptstyle g_{2}$};
\node at (8,-3.4,0) {$\scriptstyle g_{6}$};
\node at (10,-3.4,0) {$\scriptstyle g_{5}$};
\node at (12,-3.4,0) {$\scriptstyle g_{7}$};
\draw [densely dotted] (2,-3.1,0) -- (6,-3.1,0) -- (6,-2.5,0) -- (2,-2.5,0) -- cycle;
\draw [densely dotted] (8,-3.1,0) -- (10,-3.1,0) -- (10,-2.5,0) -- (8,-2.5,0) -- cycle;
\draw[->] (7,-0.5,0) -- node[right] {$\scriptstyle\upnu_{\scrJ}$}(7,-2.25,0);
\end{tikzpicture}
\end{array}\label{mut pic} 
\end{equation}
where each connected component gets reflected.  Hence, if we consider the permutation
\[
\upnu_\scrJ\upomega\colonequals
\left(  
\begin{array}{ccccccc}
1&2&3&4&5&6&7\\
\upomega(1)&\upomega(4)&\upomega(3)&\upomega(2)&\upomega(6)&\upomega(5)&\upomega(7)
\end{array}
\right)
\]
we see that $M^{\upomega}$ mutates, at the choice of summands $\scrJ$, to $M^{\upnu_\scrJ\upomega}$.  Since mutation always gives a derived equivalence, we thus obtain an equivalence
\[
\Db(\mod\End_{\scrR}(M^{\upomega}))\to
\Db(\mod\End_{\scrR}(M^{\upnu_\scrJ \upomega}))
\]
and thus, via \eqref{VdB for Type A}, an equivalence
\[
\Db(\coh \scrX^{\upomega})\to
\Db(\coh\scrX^{\upnu_\scrJ\upomega}).
\]
In this manner, we generate many equivalences between the varying $\scrX^{\upomega}$.  Three points stand out: (1) it is clear from the picture that $\upnu_{\scrJ}\upnu_{\scrJ}(M)=M$, and $\upnu_{\scrJ}(M)=M$ if and only if $\scrJ$ is componentwise symmetric, (2)  whether the above operation is a flop, or not, depends on what else gets contracted when we contract the curves in $\scrJ$; this can be checked using \ref{contraction theorem}, and (3) in the case it is a flop, the mutation functor is isomorphic to Bridgeland's flop functor, or its inverse depending on conventions.

\subsection{Back to Hyperplanes}
For all $\scrR\colonequals \mathbb{C}\llsq u,v,x,y\rrsq/(uv-f_1\hdots f_n)$, with each $f_i$  irreducible, the associated hyperplane arrangement $\scrH$ is always the Type $A$ root system $A_{n-1}$, and $\scrH^{\aff}$ its affine version \cite{HomMMP}. This is true regardless of whether $\scrR$ is isolated, or whether it admits a smooth minimal model.  As such, again whilst the combinatorics of such $\scrR$ give a reasonable intuition, they are not representative of the general case.

\begin{example}
We sketch the case $uv=f_1f_2f_3$, since the arrangement $\scrH$ is inside $\mathbb{R}^2$, which is easy to draw.  The finite arrangement $\scrH$ is the $A_2$ root system, and the minimal models of $\Spec\scrR$ are illustrated as follows.
\[
\begin{array}{ccc}
\begin{array}{c}
\begin{tikzpicture}[scale=1]
\coordinate (A1) at (135:2cm);
\coordinate (A2) at (-45:2cm);
\draw[green!60!black] (A1) -- (A2);
\draw[green!60!black] (-2,0)--(2,0);
\draw[green!60!black] (0,-2)--(0,2);
\node at (45:1.25cm) {$ \begin{tikzpicture}[xscale=0.2,yscale=0.2]
\draw[black] (-0.1,-0.04) to [bend left=25] (2.1,-0.04,0);
\draw[black] (1.9,-0.04) to [bend left=25] (4.1,-0.04,0);
\filldraw [red] (0,0) circle (1pt);
\filldraw [red] (2,0) circle (1pt);
\filldraw [red] (4,0) circle (1pt);
\node at (0,-0.7) {$\scriptstyle f_{1}$};
\node at (2,-0.7) {$\scriptstyle f_{2}$};
\node at (4,-0.7) {$\scriptstyle f_{3}$};
\end{tikzpicture}$};
\node[rotate=-50] at (115:1.25cm) {$ \begin{tikzpicture}[xscale=0.2,yscale=0.2]
\draw[black] (-0.1,-0.04) to [bend left=25] (2.1,-0.04,0);
\draw[black] (1.9,-0.04) to [bend left=25] (4.1,-0.04,0);
\filldraw [red] (0,0) circle (1pt);
\filldraw [red] (2,0) circle (1pt);
\filldraw [red] (4,0) circle (1pt);
\node at (0,-0.7) {$\scriptstyle f_{2}$};
\node at (2,-0.7) {$\scriptstyle f_{1}$};
\node at (4,-0.7) {$\scriptstyle f_{3}$};
\end{tikzpicture}$};
\node[rotate=-20] at (160:1.25cm) {$ \begin{tikzpicture}[xscale=0.2,yscale=0.2]
\draw[black] (-0.1,-0.04) to [bend left=25] (2.1,-0.04,0);
\draw[black] (1.9,-0.04) to [bend left=25] (4.1,-0.04,0);
\filldraw [red] (0,0) circle (1pt);
\filldraw [red] (2,0) circle (1pt);
\filldraw [red] (4,0) circle (1pt);
\node at (0,-0.7) {$\scriptstyle f_{2}$};
\node at (2,-0.7) {$\scriptstyle f_{3}$};
\node at (4,-0.7) {$\scriptstyle f_{1}$};
\end{tikzpicture}$};
\node at (225:1.25cm) {$ \begin{tikzpicture}[xscale=0.2,yscale=0.2]
\draw[black] (-0.1,-0.04) to [bend left=25] (2.1,-0.04,0);
\draw[black] (1.9,-0.04) to [bend left=25] (4.1,-0.04,0);
\filldraw [red] (0,0) circle (1pt);
\filldraw [red] (2,0) circle (1pt);
\filldraw [red] (4,0) circle (1pt);
\node at (0,-0.7) {$\scriptstyle f_{3}$};
\node at (2,-0.7) {$\scriptstyle f_{2}$};
\node at (4,-0.7) {$\scriptstyle f_{1}$};
\end{tikzpicture}$};
\node[rotate=-50] at (-65:1.25cm) {$ \begin{tikzpicture}[xscale=0.2,yscale=0.2]
\draw[black] (-0.1,-0.04) to [bend left=25] (2.1,-0.04,0);
\draw[black] (1.9,-0.04) to [bend left=25] (4.1,-0.04,0);
\filldraw [red] (0,0) circle (1pt);
\filldraw [red] (2,0) circle (1pt);
\filldraw [red] (4,0) circle (1pt);
\node at (0,-0.7) {$\scriptstyle f_{3}$};
\node at (2,-0.7) {$\scriptstyle f_{1}$};
\node at (4,-0.7) {$\scriptstyle f_{2}$};
\end{tikzpicture}$};
\node[rotate=-20] at (-20:1.25cm) {$ \begin{tikzpicture}[xscale=0.2,yscale=0.2]
\draw[black] (-0.1,-0.04) to [bend left=25] (2.1,-0.04,0);
\draw[black] (1.9,-0.04) to [bend left=25] (4.1,-0.04,0);
\filldraw [red] (0,0) circle (1pt);
\filldraw [red] (2,0) circle (1pt);
\filldraw [red] (4,0) circle (1pt);
\node at (0,-0.7) {$\scriptstyle f_{1}$};
\node at (2,-0.7) {$\scriptstyle f_{3}$};
\node at (4,-0.7) {$\scriptstyle f_{2}$};
\end{tikzpicture}$};
\end{tikzpicture}
\end{array}
\end{array}
\]
Each codimension one wall crossing is visibly a simple mutation as in \eqref{mut pic}, reflecting precisely one curve. Whether the minimal models on either side of any given wall are distinct, or not, depends on the choice of the $f_i$. 
\end{example}

\section{Autoequivalences and Stability Conditions}\label{autostab section}

Via tilting and mutation, it is possible to lift the combinatorial statements in the Auslander--McKay correspondence \ref{AM main} to the derived category level.  Most of this section works for an arbitrary $\scrR$, however since the statements are much more elegant if we restrict to the case when $\scrR$ is isolated, we often make this restriction.

\subsection{Spherical functors and all that}  In the setting of cDV singularities, the contraction algebra gives a natural candidate for the base of a spherical functor. Alas, when $\scrR$ is not isolated, whether this functor is spherical or not turns out to be delicate.

The construction of the functor is not the problem. Indeed, for an arbitrary crepant partial resolution $\scrX\to\Spec\scrR$ with $\scrR$ cDV, we can choose a subset $\scrI$ of curves above the origin.  As in \eqref{cont general} a certain factor of $\End_{\scrR}(N)$ is the contraction algebra, denoted here $\CA^{\scrI}$, associated to this choice. Being a quotient, there is a group homomorphism
\[
\End_{\scrR}(N)\twoheadrightarrow \CA^{\scrI}
\]
which in turn induces an exact functor between the corresponding module categories, in the reverse direction. In turn, this induces a functor
\[
\Db(\mod\CA^{\scrI})\to\Db(\mod\End_{\scrR}(N))\simeq\Db(\coh \scrX).
\]
The problem is whether this is spherical, and thus induces an autoequivalence on $\Db(\coh \scrX)$ via a certain functorial triangle.

When $\scrR$ is not isolated, there is a very precise answer to when this is the case, but it is very technical to state \cite[6.3]{DW4}, and not the easiest to check in practice.  In contrast, when $\scrR$ is isolated, and thus $\scrX\to\Spec\scrR$ is a flopping contraction, it turns out that the above functor \emph{is} always spherical.  Not only that, the whole set-up globalises as follows, and works regardless of whether the set of chosen curves contracts algebraically, or not.

\begin{thm}\label{perf and triangle}
Suppose that $\scrY\to\scrY_{\con}$ is a flopping contraction between $3$-fold quasi-projective varieties, where $\scrY$ has only Gorenstein terminal singularities in a neighbourhood of the flopping curves.  Pick a subset $\scrJ$ of these curves, and consider the universal sheaf $\scrE_\scrJ\in\coh \scrY$ from noncommutative deformation theory. Then the following hold.
\begin{enumerate}
\item $\scrE_{\scrJ}\in\coh\scrY$ is a perfect complex.
\item\label{perf and triangle 2}  There exists a functorial triangle
\[
\RHom_{\scrY}(\scrE_{\scrJ},-)\otimes_{\CA^{\scrJ}}^{\bf L}\scrE_{\scrJ}\to \Id\to\twistGen_{\scrJ}\to.
\]
where $\twistGen_{\scrJ}$ is an autoequivalence.
\end{enumerate}
\end{thm}

The first statement, that $\scrE_\scrJ$ is a perfect complex, is perhaps the most difficult part to prove, and it heavily relies on the fact that $\scrY$ has only Gorenstein terminal singularities in a neighbourhood of the flopping curves.  The importance of $\scrE_\scrJ$ being perfect should really be understood as offering the `hidden smoothness' in the derived category of $\scrY$ that allows us to avoid using stacks.

\medskip
The twist functors in \ref{perf and triangle}\eqref{perf and triangle 2} are all well and lovely, but the functional triangles are as good as useless if you want to know which subgroup of $\mathrm{Auteq}\Db(\coh \scrY)$ that the twists generate.  The groupoid picture explained in the next subsection fixes this problem, but it also does more.  There are two points: (1) in dimension three it is the \emph{flop functors} that braid\footnote{the length of the braid relation need not be three, however}, not the spherical twists, and (2) the groupoid picture allows us to visually see the limitations of only studying compositions of flop functors and line bundle twists.  The autoequivalence group is a \emph{lot} larger than that.

\subsection{Groupoids} 
Given an intersection arrangement $\scrH$ or $\scrH^{\aff}$ from \S\ref{int arr section},  there is an associated Deligne groupoid, defined combinatorially as follows.  For each chamber there is a corresponding vertex, and two vertices are connected by an arrow if and only if the chambers are adjacent.
\[
\includegraphics[angle=0,scale = 0.7]{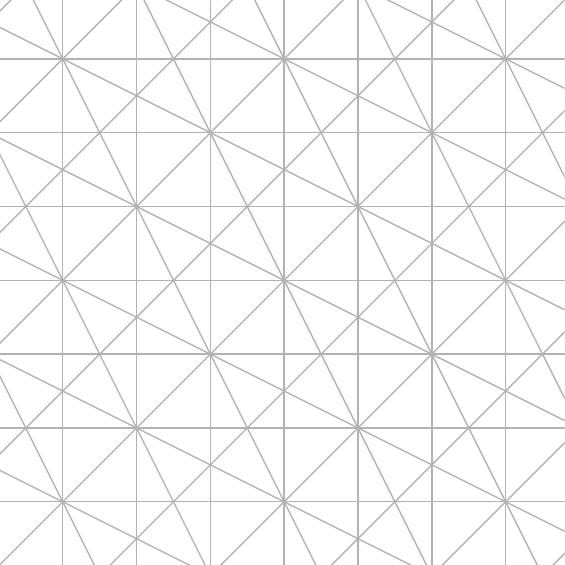}
\qquad
\includegraphics[angle=0,scale = 0.7]{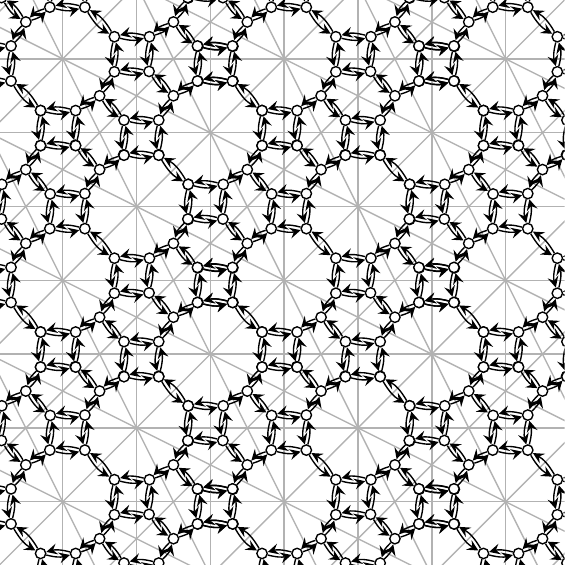}
\]
We then add relations given by identifying shortest paths, then formally invert all arrows. This gives rise to the Deligne groupoid, written $\mathds{G}_{\scrH}$ in the finite case, or $\mathds{G}_{\scrH^{\aff}}$ in the infinite case illustrated above.

\medskip
Restricting to the case $\scrR$ is isolated cDV in order to obtain the most elegant statement, there is another way to build a groupoid.  As in \S\ref{AM section} by choosing a minimal model we obtain arrangements $\scrH$ and $\scrH^{\aff}$. Since by \ref{bijection main isol} there is a bijection between MM and chambers, and furthermore wall crossing corresponds to mutation, we can build groupoids $\scrG_{\scrH}$ and $\scrG_{\scrH^{\aff}}$, generated by the following:
\begin{itemize}
\item Each chamber has associated $M$, thus $\Db(\mod\End_{\scrR}(M))$.  In the case $\scrG_{\scrH}$ this is indexed over $M\in(\mathrm{MM}\,\scrR)\cap(\CM\scrR)$, for $\scrG_{\scrH^{\aff}}$ it is indexed over $M\in\mathrm{MM}\,\scrR$.
\item Each arrow has an associated \emph{mutation} equivalence from \cite{IW3}.
\end{itemize}
The following asserts that the mutation functors satisfy the shortest path relations.

\begin{thm}[\cite{HW1, IW9}]\label{aff groupoid}
When $\scrR$ is isolated \textnormal{cDV}, there exists functors $\mathds{G}_{\scrH}\to\scrG_{\scrH}$ and $\mathds{G}_{\scrH^{\aff}}\to\scrG_{\scrH^{\aff}}$.  The first is faithful.
\end{thm}

The second functor is conjecturally faithful, but this is only known in very specific Type $A$ cases.  However, it is really the corollaries of \ref{aff groupoid} that interest us the most.  The first is really just a special case of \ref{aff groupoid}, since some of the mutation functors correspond to the flop functors, via \cite[4.2]{HomMMP}.

\begin{cor}[Flops Braid \cite{DW3}]
Suppose that $\scrX\to \Spec\scrR$ is a flopping contraction where $\scrR$ is isolated \textnormal{cDV}. Then for any curves $\mathrm{C}_1$ and $\mathrm{C}_2$ above the origin, 
\[
\underbrace{\mathsf{Flop}_1\circ\mathsf{Flop}_2\circ\mathsf{Flop}_1\circ\cdots}_{d}
\cong
\underbrace{\mathsf{Flop}_2\circ\mathsf{Flop}_1\circ\mathsf{Flop}_2\circ\cdots}_{d}
\]
where $d$ equals either $2, 3, 4, 5, 6$ or $8$.  All such $d$ can be realised.
\end{cor}

The second corollary is immediate from \ref{aff groupoid} by passing to vertex groups of the associated groupoids.  As is well known, the vertex group of the Deligne groupoid is isomorphic to the fundamental groups of the corresponding complexified complement; as in the introduction denote these $\scrZ$ and $\scrZ_{\aff}$ respectively.

\begin{cor}\label{actions on Db}
Suppose that $\scrX\to\Spec\scrR$ is a flopping contraction, with $\scrR$ isolated \textnormal{cDV}.  Then there are group homomorphisms
\[
\begin{tikzpicture}
\node (A1) at (0,0) {$\uppi_1(\scrZ)$};
\node (A2) at (0,-1.5) {$\uppi_1(\scrZ_{\aff})$};
\node (B) at (2.5,0) {$\mathrm{Auteq}\Db(\coh \scrX).$};
\draw[->] (A1) -- node[above] {$\scriptstyle\upvarphi$}(B);
\draw[->] (A2) --  node[gap]  {$\scriptstyle\widetilde{\upvarphi}$}(B);
\draw[->] (A1) -- (A2);
\end{tikzpicture}
\]
\end{cor}

\subsection{Stability}
Continuing with the assumption that $\scrX\to\Spec\scrR$ is a flopping contraction, with $\scrR$ isolated \textnormal{cDV}, consider the following two subcategories of $\Db(\coh \scrX)$.
\begin{align*}
\scrC&= \{ \scrF\in\Db(\coh \scrX)\mid \mathbf{R
}f_*\scrF=0\} \\
\scrD&= \{ \scrF\in\Db(\coh \scrX)\mid \Supp\scrF\subseteq \mathrm{C} \}.
\end{align*}
The category $\scrD$ forms a local model for a projective CY $3$-fold.  Following work of Bridgeland for Kleinian singularities \cite{B3}, we should view $\scrC$ as the `finite type' category, and $\scrD$ as the affine version.  

We do not give an overview of stability conditions here, except to remark that the full spaces of stability conditions $\mathrm{Stab}\,\scrC$ and $\mathrm{Stab}\,\scrD$ are complex manifolds.  So as to match the choice of level in \eqref{level}, below we \emph{normalise} stability conditions on $\scrD$ with respect to the scaling action, and so instead consider $\mathrm{Stab}_n\scrD$.

\begin{thm}[\cite{HW}]\label{main stab}
Suppose $\scrX\to\Spec\scrR$ is a crepant partial resolution, with $\scrR$ isolated \textnormal{cDV}, and associate finite $\scrH$ and infinite $\scrH^{\mathsf{aff}}$ by slicing.  Then there are connected components $\mathrm{Stab}^\circ\scrC$ and $\mathrm{Stab}^\circ_n\scrD$ such that the  forgetful maps
\begin{align*}
\mathrm{Stab}^\circ\scrC&\to\mathbb{C}^n\backslash\scrH_{\mathbb{C}}\\
\mathrm{Stab}^\circ_n\scrD&\to\mathbb{C}^n\backslash(\scrH^{\aff})_{\mathbb{C}}
\end{align*}
are regular covering maps, with Galois groups $\uppi_1(\scrZ)$ and $\uppi_1(\scrZ_{\aff})$ respectively.   The first is universal.
\end{thm}

There are many consequences of \ref{main stab}, not least a computation of the string K\"ahler moduli space for $3$-fold flopping contractions; see \cite[\S7]{HW} for full details.

\begin{example}
In the following example, of a certain two-curve $E_8$ flop $\scrX\to\Spec\scrR$, inside $\scrH^{\aff}$ the movable cone (containing 16 chambers) is highlighted in green, as are all its translations by $\mathrm{Pic}(\scrX)\cong\mathbb{Z}^{\oplus 2}$.
\[
\includegraphics[angle=0,scale = 0.4]{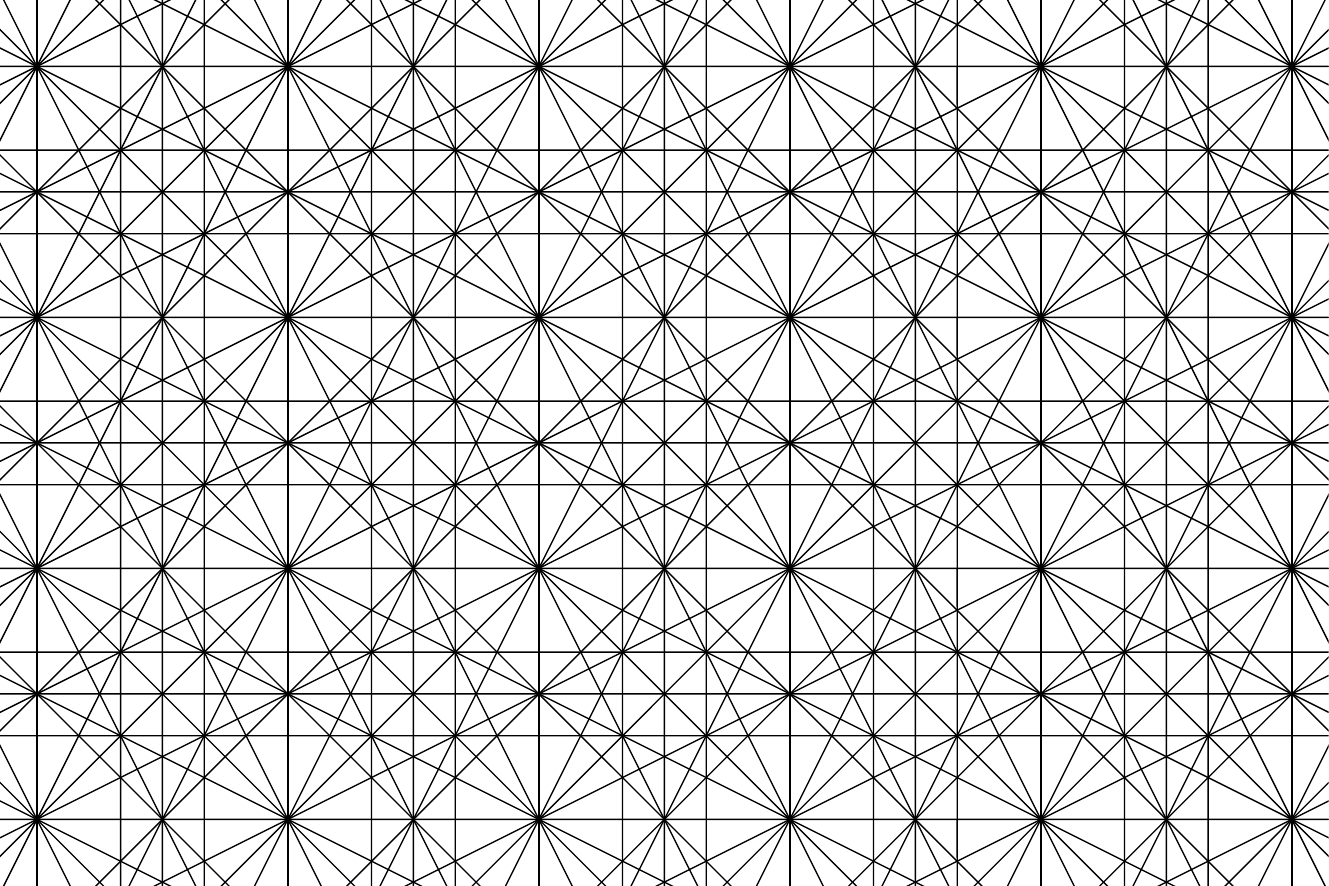}
\quad
\includegraphics[angle=0,scale = 0.4]{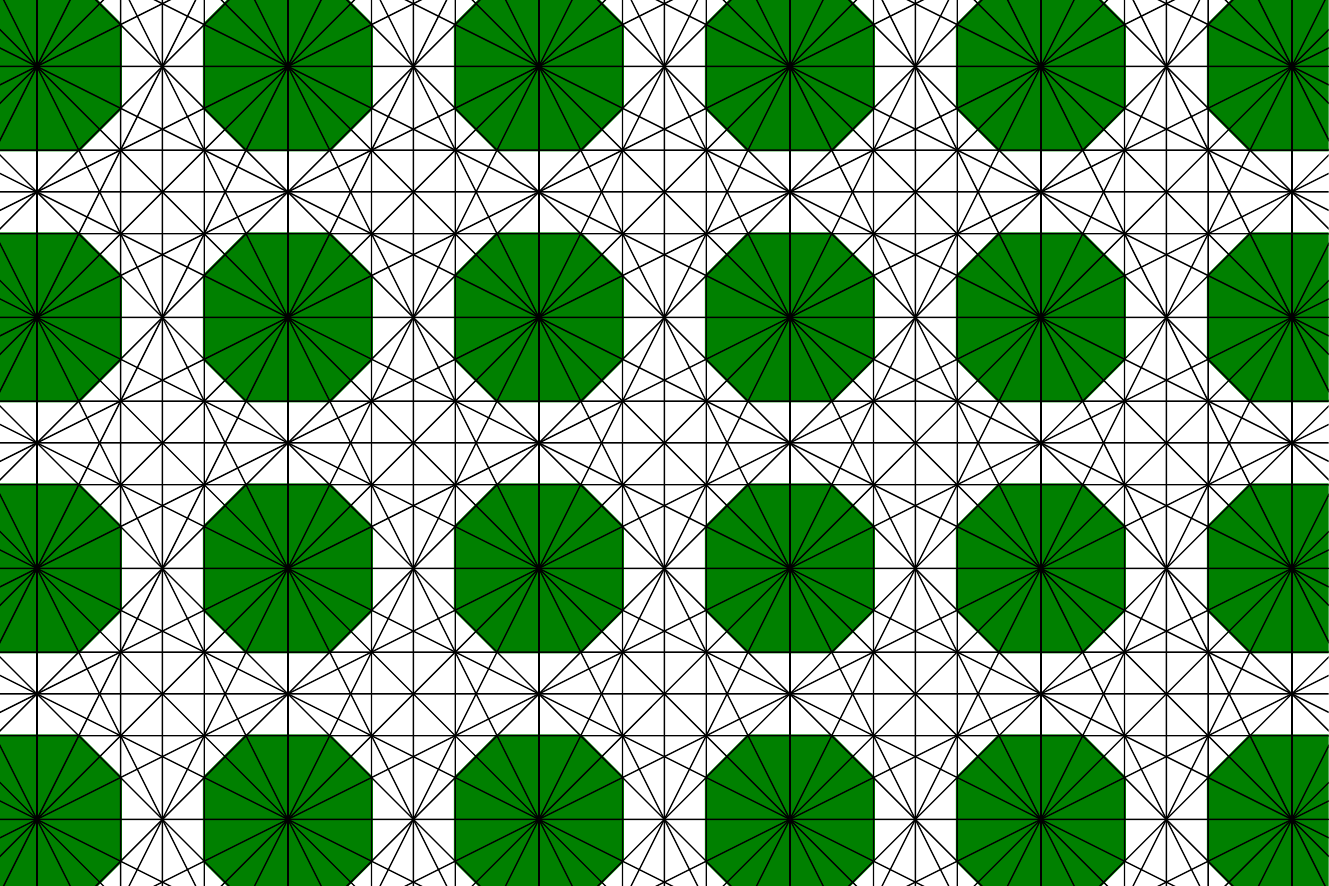}
\]
The areas between the green islands do not have an obvious birational geometry interpretation, but monodromy around hyperplanes in those areas still gives rise to elements in  $\mathrm{Image}(\widetilde{\upvarphi})$, where $\widetilde{\upvarphi}$ is from \ref{actions on Db}.  This demonstrates, visually, why the group consisting of compositions of flop functors and line bundle twists is in general much smaller than $\mathrm{Image}(\widetilde{\upvarphi})\rtimes\mathrm{Pic}\,\scrX$, and thus in particular is much smaller than $\mathrm{Auteq}\Db(\coh \scrX)$.  

This phenomena is only visible once we pass to the more complicated examples of flops: indeed, when $\scrX\to\Spec\scrR$ slices to a minimal resolution,   compositions of flop functors and line bundle twists  suffice \cite{TodaResPub}.
\end{example}


\begin{thebibliography}{DW4}

\bibitem[A1]{AugustThesis}
J.~August, \emph{Tilting theory of contraction algebras}, PhD thesis, University of Edinburgh, 2019, \href{https://era.ed.ac.uk/handle/1842/35985}{link}.

\bibitem[A2]{August1}
J.~August, \emph{On the finiteness of the derived equivalence classes of some stable endomorphism rings}, Math.\ Z.\ \textbf{296} (2020), no.\ 3-4, 1157--1183. 

\bibitem[A3]{August2}
J.~August, \emph{The tilting theory of contraction algebras}, Adv.\ Math.\ \textbf{374} (2020), 107372, 56 pp.

\bibitem[A4]{Aus2}
M.~Auslander, \emph{Isolated singularities and existence of almost split sequences}. Representation theory, II (Ottawa, Ont., 1984), 194--242, Lecture Notes in Math., 1178, Springer, Berlin, 1986. 

\bibitem[A5]{Aus3}
M.~Auslander, \emph{Rational singularities and almost split sequences}. Trans.\ Amer.\ Math.\ Soc.\ \textbf{293} (1986), no.~2, 511--531.


\bibitem[B1]{Booth}
M.~Booth, \emph{Singularity categories via the derived quotient}, Adv.\ Math.\ \textbf{381} (2021), 107631.

\bibitem[B2]{Bridgeland}
T.~Bridgeland, \emph{Flops and derived categories}. Invent.\ Math.\ \textbf{147} (2002), no.~3, 613--632.


\bibitem[B3]{B3} T.~Bridgeland, {\it Stability conditions and Kleinian singularities},  Int.\ Math.\ Res.\ Not.\ IMRN (2009), no.~21, 4142--4157.

\bibitem[BW]{BW} G.~Brown and M.~Wemyss, \emph{Gopakumar--Vafa invariants do not determine flops}, Comm.\ Math.\ Phys.\ \textbf{361} (2018), no.~1, 143--154.

\bibitem[BW2]{BW2} G.~Brown and M.~Wemyss, \emph{Local Normal Forms of Noncommutative Functions}, in preparation.

\bibitem[BKL]{BKL}
J.~Bryan, S.~Katz and N.~Leung, \emph{Multiple covers and the integrality conjecture for rational curves in Calabi-Yau threefolds}, J.\ Algebraic Geom.\ \textbf{10} (2001), no.~3, 549--568. 

\bibitem[BIKR]{BIKR}
I.~Burban, O.~Iyama, B.~Keller and I.~Reiten, \emph{Cluster tilting for one-dimensional hypersurface singularities}, Adv.\ Math.\ \textbf{217} (2008), no.~6, 2443--2484.


\bibitem[DH]{DH}
H.~Dao and C.~Huneke, \emph{Vanishing of Ext, cluster tilting modules and finite global dimension of endomorphism rings}, Amer.\ J.\ Math.\ \textbf{135} (2013), no.~2, 561--578.

\bibitem[D]{Davison}
B.~Davison, \emph{Refined invariants of finite-dimensional Jacobi algebras}. 
\href{https://arxiv.org/abs/1903.00659}{\sf arXiv:1903.00659}.

\bibitem[DW1]{DW1}
W.~Donovan and M.~Wemyss, \emph{Noncommutative deformations and flops}, Duke Math.\ J.\ \textbf{165} (2016), no.~8, 1397--1474. 

\bibitem[DW2]{DW2}
W.~Donovan and M.~Wemyss, \emph{Contractions and deformations}, Amer.\ J.\ Math.\ \textbf{141} (2019), no.~3, 563--592.

\bibitem[DW3]{DW3}
W.~Donovan and M.~Wemyss, \emph{Twists and braids for general 3-fold flops}, J.\ Eur.\ Math.\ Soc.\ (JEMS), \textbf{21} (2019), no.~6, 1641--1701.

\bibitem[DW4]{DW4}
W.~Donovan and M.~Wemyss, \emph{Noncommutative enhancements of contractions}, Adv.\ Math.\ \textbf{344} (2019), 99--136.

\bibitem[E1]{Eriksen2}
E.~Eriksen, \emph{An introduction to noncommutative deformations of modules}, Noncommutative algebra and geometry, 90--125, Lect.\ Notes Pure Appl.\ Math., \textbf{243}, Chapman and Hall/CRC, Boca Raton, FL, 2006.

\bibitem[E2]{Eriksen}
E.~Eriksen, \emph{Computing noncommutative deformations of presheaves and sheaves of modules}, Canad.\ J.\ Math.\ \textbf{62} (2010), no.~3, 520--542.

\bibitem[HW1]{HW1} Y.~Hirano and M.~Wemyss, {\it Faithful actions from hyperplane arrangements}, Geom.\  Topol.\ \textbf{22} (2018), no.~6, 3395--3433. 

\bibitem[HW2]{HW} 
Y.~Hirano and M.~Wemyss, {\it Stability conditions for $3$-fold flops}, \href{https://arxiv.org/abs/1907.09742}{\sf arXiv:1907.09742}.

\bibitem[H]{Hua}
Z.~Hua, \emph{Contraction algebra and singularity of three-dimensional flopping contraction}, Math.\ Z.\ \textbf{290} (2018), no.~1-2, 431--443.

\bibitem[HK]{HuaKeller}
Z.~Hua and B.~Keller, \emph{Cluster categories and rational curves}, \href{https://arxiv.org/abs/1810.00749}{\sf arXiv:1810.00749}.

\bibitem[HT]{HuaToda}
Z.~Hua and Y.~Toda, \emph{Contraction algebra and invariants of singularities}, Int.\ Math.\ Res.\ Not.\ IMRN 2018, no.~10, 3173--3198.

\bibitem[IW1]{IW3}
O.~Iyama and M.~Wemyss, \emph{Maximal modifications and Auslander--Reiten duality for non-isolated singularities}, Invent.\ Math. \textbf{197} (2014), no.~3, 521--586. 

\bibitem[IW2]{IW5}
O.~Iyama and M.~Wemyss, \emph{Singular derived categories of $\mathbb{Q}$-factorial terminalizations and maximal modification algebras}, Adv.\ Math.\ \textbf{261} (2014), 85--121.

\bibitem[IW3]{IW6}
O. Iyama and M. Wemyss, \emph{Reduction of triangulated categories and Maximal Modification Algebras for $cA_n$ singularities}, J.\ Reine Angew.\ Math.\ \textbf{738} (2018), 149--202.

 \bibitem[IW4]{IW9}
O.~Iyama and M.~Wemyss, \emph{Tits cones intersections and applications}, \href{https://www.maths.gla.ac.uk/~mwemyss/MainFile_for_web.pdf}{\sf preprint}. 

\bibitem[I]{Iyudu}
N.~Iyudu, \emph{Classification of contraction algebras and pre-Lie algebras associated to braces and trusses}, \href{https://arxiv.org/abs/2008.06033}{\sf arXiv:2008.06033}.

\bibitem[IS1]{IS}
N.~Iyudu and S.~Shkarin, \emph{Potential algebras with few generators}, Glasg.\ Math.\ J.\ \textbf{62} (2020), no.~S1, S28--S76.

\bibitem[IS2]{IS2}
N.~Iyudu and A.~Smoktunowicz, \emph{Golod-Shafarevich-type theorems and potential algebras}, Int.\ Math.\ Res.\ Not.\ IMRN 2019, no.~15, 4822--4844.

\bibitem[K1]{Katz}
S.~Katz, \emph{Genus zero Gopakumar-Vafa invariants of contractible curves}, J.\ Differential Geom.\ \textbf{79} (2008), no.~2, 185--195.


\bibitem[K2]{KawamataCone}
 Y.~Kawamata, \emph{On the cone of divisors of Calabi-Yau fiber spaces}, Internat.\ J.\ Math.\ \textbf{8} (1997), no.~5, 665--687.
 
\bibitem[K3]{Kawamata}
Y.~Kawamata, \emph{On multi-pointed non-commutative deformations and Calabi-Yau threefolds}, Compos.\ Math.\ \textbf{154} (2018), no.~9, 1815--1842.

\bibitem[K4]{Kawamata2}
Y.~Kawamata, \emph{Non-commutative deformations of simple objects in a category of perverse coherent sheaves}, Selecta Math.\ (N.S.) \textbf{26} (2020), no.~3, Paper No.~43, 22 pp. 

\bibitem[K5]{Kawamata3}
Y.~Kawamata, \emph{Non-commutative deformations of perverse coherent sheaves and rational curves}, \href{https://arxiv.org/abs/2006.09547}{\sf arXiv:2006.09547}.

\bibitem[K6]{KollarFlops}
J.~Koll\'ar, \emph{Flops}. Nagoya Math.\ J.\ \textbf{113} (1989), 15--36.

\bibitem[KM]{KollarMori}
J.~Koll\'ar and S.~Mori, \emph{Classification of three-dimensional flips}, J.\ Amer.\ Math.\ Soc.\ \textbf{5} (1992), no.~3, 533--703.

\bibitem[L1]{Laudal}
O.~A.~Laudal, \emph{Noncommutative deformations of modules}. The Roos Festschrift volume, 2.\ Homology Homotopy Appl.\ \textbf{4} (2002), no.~2, part 2, 357--396.

\bibitem[L2]{LeuschkeSurvey}
G.~J.~Leuschke, \emph{Non-commutative crepant resolutions: scenes from categorical geometry}, Progress in commutative algebra \textbf{1}, 293--361, de Gruyter, Berlin, 2012. 

\bibitem[MW]{MakWu}
C.~Y.~Mak and W.~Wu, \emph{Dehn twists and Lagrangian spherical manifolds}, Selecta Math.\ (N.S.) \textbf{25} (2019), no.~5, Paper No.~68, 85 pp.

\bibitem[M1]{Markushevich}
D.~G.~Markushevich, \emph{Canonical singularities of three-dimensional hypersurfaces},(Russian) Izv.\ Akad.\ Nauk SSSR Ser.\ Mat.\ \textbf{49} (1985), no.~2, 334--368, 462.

\bibitem[M2]{McKay}
J.~McKay, \emph{Graphs, singularities, and finite groups}, Proc.\ Sympos.\ Pure  Math.\ \textbf{37} (1980), 183--186.

\bibitem[N]{Nag}
K.~Nagao, \emph{Derived categories of small toric Calabi-Yau 3-folds and curve counting invariants}, Q.\ J.\ Math.\ \textbf{63} (2012), no.~4, 965--1007.

\bibitem[P]{Pinkham}
H.~Pinkham, \emph{Factorization of birational maps in dimension 3}, Singularities (P. Orlik, ed.), Proc.\ Symp.\ Pure Math., vol.\ 40, Part 2, 343--371, American Mathematical Society, Providence, 1983.

\bibitem[R1]{ReidCanonical}
M.~Reid, \emph{Canonical 3-folds}, Journ\'{e}es de G\'{e}ometrie Alg\'{e}brique d'Angers, Juillet 1979/Algebraic Geometry, Angers, 1979, pp. 273--310, Sijthoff \& Noordhoff, Alphen aan den Rijn--Germantown, Md., 1980. 

\bibitem[R2]{Pagoda}
M.~Reid, \emph{Minimal models of canonical 3-folds}, Algebraic varieties and analytic varieties (Tokyo, 1981), 131--180, Adv. Stud. Pure Math., 1, North-Holland, Amsterdam, 1983.

\bibitem[R3]{YPG}
M.~Reid, \emph{Young person's guide to canonical singularities}, Algebraic geometry, Bowdoin, 1985 (Brunswick, Maine, 1985), 345--414, Proc.\ Sympos.\ Pure Math., \textbf{46}, Part 1, Amer.\ Math.\ Soc., Providence, RI, 1987. 

\bibitem[ST]{ST}
P.~Seidel and R.~Thomas, \emph{Braid group actions on derived categories of coherent sheaves}, Duke Math.\ J.\ \textbf{108} (2001), no.~1, 37--108.

\bibitem[SW]{SW}
I.~Smith and M.~Wemyss, \emph{Double bubble plumbings and two-curve flops}, \href{https://arxiv.org/abs/2010.10114}{\sf arXiv:2010.10114}.


\bibitem[T1]{TodaResPub}
Y.~Toda, \emph{Stability conditions and crepant small resolutions}, Trans.\ Amer.\ Math.\ Soc.\ \textbf{360} (2008), no.~11, 6149--6178.

\bibitem[T2]{TodaGV}
Y.~Toda, \emph{Non-commutative width and Gopakumar-Vafa invariants},
Manuscripta Math.\ \textbf{148} (2015), no.~3-4, 521--533. 

\bibitem[T3]{TodaUtah}
Y.~Toda, \emph{Non-commutative deformations and Donaldson-Thomas invariants}, Algebraic geometry: Salt Lake City 2015, 611--631, Proc.\ Sympos.\ Pure Math., 97.1, Amer.\ Math.\ Soc., Providence, RI, 2018.

\bibitem[V1]{Viehweg}
E.~Viehweg, \emph{Rational singularities of higher dimensional schemes},
Proc.\ Amer.\ Math.\ Soc.\ \textbf{63} (1977), no.~1, 6--8. 

\bibitem[V2]{VdB1d}
M.~Van den Bergh, \emph{Three-dimensional flops and noncommutative rings}, Duke Math.\ J.\ \textbf{122} (2004), no.~3, 423--455. 

\bibitem[V3]{VdBNCCR}
M.~Van den Bergh, \emph{Non-commutative crepant resolutions}, The legacy of Niels Henrik Abel, 749--770, Springer, Berlin, 2004.

\bibitem[vG]{Okke}
O.~van Garderen, \emph{Donaldson-Thomas invariants of length 2 flops}, \href{https://arxiv.org/abs/2008.02591}{\sf arXiv:2008.02591}.

\bibitem[W1]{WemyssSurvey}
M.~Wemyss, \emph{Noncommutative resolutions}, Noncommutative algebraic geometry, 23--306, Math.\ Sci.\ Res.\ Inst.\ Publ., \textbf{64}, Cambridge Univ.\ Press, New York, 2016. 

\bibitem[W2]{HomMMP}
M.~Wemyss, \emph{Flops and Clusters in the Homological Minimal Model Program}, Invent.\ Math.\ \textbf{211} (2018), no.~2, 435--521.


\end{thebibliography}
\end{document}